\documentclass[twocolumn]{autart}    

\usepackage{enumitem}

\usepackage{graphicx}          
\usepackage{graphics}
\usepackage{subfigure}

\usepackage{amsmath,amssymb,amsfonts}

\allowdisplaybreaks
\usepackage{varioref}
\usepackage{hyperref}

\usepackage{color}

\edef\endfrontmatter{ 				
  \unexpanded\expandafter{\endfrontmatter}
  \noexpand\endNoHyper 
}

\usepackage{natbib}        
\newcommand{\lmm}{}
\newcommand{\kl}{}
\newcommand{\ml}{}
\newcommand{\icl}{}
\newcommand{\il}{}
\newcommand{\ill}{}
\newcommand{\nkl}{}
\newcommand{\mm}{}
\newcommand{\mml}{}
\newcommand{\li}{}
\newcommand{\ic}{}

\newcommand{\hiro}{}
\newcommand{\hirot}{}
\newcommand{\irr}{}
\newcommand{\irrr}{}
\newcommand{\ir}{}
\newcommand{\lir}{}
\newcommand{\hu}{}
\newcommand{\lr}{}

\newcommand{\lrrr}{}

\usepackage[prependcaption,colorinlistoftodos]{todonotes}

\newcommand{\todoiny}[1]{}
\newcommand{\todoing}[1]{}
\newcommand{\todoKL}[1]{}
\newcommand{\todoinc}[1]{}
\newcommand{\todoino}[1]{}
\newcommand{\revise}{}
\newcommand{\icr}{}
\newcommand{\icrr}{}
\newcommand{\iclr}{}
\newcommand{\icc}{}
\newcommand{\mmli}{}

\newtheorem{theorem}{Theorem}
\newtheorem{lemma}{Lemma}

\newtheorem{remark}{Remark}

\newtheorem{proposition}{Proposition}

\begin{document}

\begin{frontmatter}

\title{Convergence Rate Bounds for the Mirror Descent Method:\\ IQCs,  Popov Criterion and Bregman Divergence
\thanksref{footnoteinfo}
} 

\thanks[footnoteinfo]{A preliminary version of this paper was presented in the \ic{$61$st} IEEE Conference on Decision and Control (CDC 2022).
}

\author[Cam,Hiroshima]{Mengmou Li},
\author[Cam,Leicester]{Khaled~Laib},
\author[Titech]{Takeshi~Hatanaka},
\author[Cam]{Ioannis~Lestas}

\address[Cam]{Department of Engineering, University of Cambridge, Trumpington Street, Cambridge CB2 1PZ, United Kingdom (email: mmli.research@gmail.com; kl507@cam.ac.uk; icl20@cam.ac.uk).}        
\address[Hiroshima]{\hiro{Graduate School of Advanced Science and Engineering, Hiroshima University, 1-2-1 Kagamiyama, Higashi-Hiroshima, Japan (email: mengmou@hiroshima-u.ac.jp).}}
\address[Titech]{Department of Systems and Control Engineering, School of Engineering, Tokyo Institute of Technology, S5-16, 2-12-1 Ookayama, Meguro-ku, Tokyo, Japan (email: hatanaka@sc.e.titech.ac.jp).} 
\address[Leicester]{\hiro{School of Engineering, University of Leicester, Leicester, LE1 7RH, United Kingdom (email: kl314@le.ac.uk).}}             


\begin{abstract}                          
This paper \lmm{presents a comprehensive} convergence analysis \ill{for the mirror descent (MD) method, a \lmm{widely used} algorithm in convex optimization.} 
\ic{The key feature of this algorithm is that it provides a generalization \li{of} classical gradient-based methods 
via the use of generalized distance-like functions, which are formulated using the Bregman divergence. Establishing convergence rate bounds for this algorithm is in general a non-trivial problem due to the lack of monotonicity properties in the composite nonlinearities involved.
In \icc{this} paper, we show that the 
Bregman divergence from the optimal solution,  
which is commonly used as a Lyapunov function for this algorithm, is a special case of Lyapunov functions that follow when the Popov criterion is applied to an appropriate reformulation of the \revise{MD dynamics}. This is then used
as a basis 
to construct 
an integral quadratic constraint (IQC) framework through which convergence rate bounds with reduced conservatism can be deduced.}
\ml{We also \il{illustrate via examples that the 
convergence rate \ic{bounds derived 
can be} tight.}}
\todoing{\st{Last sentence not very clear}}
\end{abstract}

\end{frontmatter}

\section{Introduction}
The mirror descent (MD) method was initially proposed by \cite{nemirovskij1983problem} for solving constrained convex optimization problems. By choosing a Bregman distance function in place of the Euclidean distance to reflect the geometry of the constraint sets, it generalizes the gradient descent (GD) method from the Euclidean space to Hilbert and Banach \ic{spaces (\cite{bubeck2015convex}).}
\mm{
Attempts have been made to unify or incorporate the MD method with various algorithms, such as the nonlinear projected subgradient algorithms (\cite{beck2003mirror}), online learning algorithms (\cite{raskutti2015information}), accelerated algorithms (\cite{krichene2015accelerated,wibisono2016variational}), and dual averaging (\cite{juditsky2022unifying}). Due to its applicability in non-Euclidean space, it has received considerable research attention in many contexts, such as stochastic optimization (\cite{duchi2012ergodic,nedic2014stochastic,borovykh2022stochastic}), distributed optimization (\cite{yuan2018optimal,doan2018convergence,sun2022centralized}),  and machine learning (\cite{mertikopoulos2018optimistic}).
}

\lmm{
In contrast to the abundant literature on well-known algorithms such as gradient descent and Nesterov's accelerated algorithms, there is a lack of references characterizing a tight bound on the convergence rate for the MD method with a constant stepsize. }
\revise{
While \cite{lu2018relatively} adopt the \icrr{notion of} relative smoothness/convexity and \cite{sun2022centralized} use Lyapunov functions \icrr{formulated via} quadratic constraints \revise{(QCs)}, respectively, to obtain convergence rate bounds for \iclr{the MD method} under a constant stepsize, these bounds are not tight in general. }
It is important to address \icc{these problems}, as diminishing stepsizes usually lead to slower convergence than constant stepsizes, and {\icrr{tight convergence rate bounds are \mmli{crucial} in the analysis and applications of optimization algorithms.}

\ill{\ic{When bounds} on the convergence rate need to be established, it is important to have systematic methods that allow to construct Lyapunov functions with more advanced structures, or allow via other means to deduce convergence rates with reduced conservatism. It has been pointed out in the optimization literature that integral quadratic constraint (IQC) framework (\cite{megretski1997system}) can be a useful tool in this direction 
(\revise{\cite{lessard2016analysis,dhingra2018proximal,scherer2023optimization}}).
However, \revise{there has been no application of IQCs} to the case of MD algorithm, as the MD dynamics involve the composition of two nonlinearities that correspond to monotone operators, with this composition not preserving these monotonicity properties.}

\ml{In the analysis of MD algorithm, the \ill{commonly used} Lyapunov function is the Bregman divergence function representing a generalized distance between the decision variable and the optimal solution.} The Bregman divergence function was introduced by \cite{bregman1967relaxation} to find the intersection of convex sets. It has \icrr{had} applications in the analysis of distributed optimization (\cite{li2020input}), port-Hamiltonian systems (\cite{jayawardhana2007passivity}), \ill{equilibrium independent stability (\cite{simpson2018hill}), and power systems (\cite{de2017bregman,monshizadeh2019secant}), in addition to the analysis of the MD algorithm.

\icrr{Despite the nice geometric interpretation of the Bregman divergence in the context of the MD method, no link has been established thus far
between the use of the Bregman divergence from the optimal solution as a Lyapunov function and
robustness analysis tools. Such connections are important as they can motivate wider classes of multipliers through which less conservative convergence rate bounds can be established.}

\icrr{One of our main contributions in this manuscript is to show that there is indeed such an inherent connection, whereby  the \iclr{Bregman-based Lyapunov function for the MD method is a
special case of Lyapunov functions that follow when the Popov criterion is appropriately \iclr{applied.}}
The problem \iclr{formulation} that leads to this result can therefore be used as a basis to develop an IQC framework for establishing improved rate bounds for the MD method.}

\ill{
Our contributions in this paper can be summarized as follows:
\begin{enumerate}[label=\roman*)]

	\item \icrr{We \ic{\icrr{show} 
that the use 
of the} Bregman divergence from the optimal solution as a Lyapunov function for the \ic{continuous-time} MD method is a special case of Lyapunov functions \ic{that follow} 
from the Popov criterion, when this is applied to an appropriate reformulation of the \revise{\icrr{MD dynamics}}.
}

\item
\icrr{
\revise{\icrr{We use the reformulation of the MD dynamics that leads to the result above
to develop \mmli{an} IQC framework 
for deriving convergence rate bounds.}}
We employ conic combinations of Popov IQCs and other types of IQCs 
\icrr{to derive} convergence rate bounds for the discrete-time MD method. 
\ic{\revise{These bounds are \icrr{less conservative than the} current state of the art and their tightness} is also illustrated \icrr{via} examples.}} 

\item \ic{We also extend our results to the case of the \revise{projected} MD \icrr{algorithm, where projections are additionally included in the dynamics.
 \revise{\icrr{A 
 framework is provided for establishing convergence rate bounds for the projected MD method
 with a constant stepsize.}}
 In particular, the} convergence rate is analyzed via additional 
IQCs associated with the projection operator and the repeated nonlinearities formulated.}

\end{enumerate}
}
\ill{The convergence rate bounds deduced are formulated as solutions to linear matrix inequalities (LMIs) in both discrete and continuous time.}
\mm{Compared to the preliminary version of this paper \cite{li2022convergence}, we \revise{obtain
\icrr{an analytical expression for the
largest convergence rate bound} in continuous time, and}
include the complete proofs of the convergence analysis \ic{for} the unconstrained MD algorithm. Moreover, we \icrr{analyze 
the convergence rate of the \revise{projected} MD algorithm.}} 

The rest of this paper is organized as follows.
In Section~\ref{Preliminaries}, \icl{preliminaries on} the MD method and \icl{IQCs are provided.} The continuous-time and discrete-time MD methods are analyzed via IQCs in Section~\ref{Continuous-time mirror descent method} and Section~\ref{Discrete-time mirror descent method}, respectively.
\revise{In addition, the projected variant of the discrete-time MD algorithm is \icl{analysed} in Section~\ref{Projected MD and convergence rate}.}
In Section~\ref{Numerical Examples}, numerical examples are given to verify our results.
Finally, the paper is concluded in Section~\ref{Conclusion}.



\section{Preliminaries}\label{Preliminaries}
\subsection{Notation}
Let $\mathbb{R}$, $\mathbb{Z}$, and $\mathbb{Z}_{+}$ denote the set of real numbers, integers, and nonnegative integers, respectively. Let $I_d$ and $0_d$ denote the $d \times d$ identity matrix and zero matrix, respectively. Their subscripts can be omitted if it is clear from the context. \ml{The notation $\textup{diag} (\alpha_1, \ldots, \alpha_d)$ denotes a $d \times d$ diagonal matrix with $\alpha_i$ on its $i$-th diagonal entry.}
Let $\mathbf{RH}_{\infty}$ be the set of proper real rational functions without poles in the closed right-half plane. The set of $m \times n$ matrices with elements in $\mathbf{RH}_{\infty}$ is denoted $\mathbf{RH}_{\infty}^{m \times n}$. Let $\mathbf{L}_{2}^{m} [0, \infty)$ be the Hilbert space of all square integrable and Lebesgue measurable functions $f: [0, \infty) \rightarrow \mathbb{R}^{m}$. It is a subspace of $\mathbf{L}_{2e}^{m}[0, \infty)$ whose elements only need to be integrable on finite intervals. Let ${l}_{2}^{m}(\mathbb{Z}_{+})$ be the set of all square summable sequences $f : \mathbb{Z}_{+} \rightarrow \mathbb{R}^{m}$.
\ml{Given a Hermitian matrix $H (j\omega)$, $H^*(j\omega) : = H^T(-j\omega)$ represents its conjugate transpose and \hiro{$\left\{ H (j\omega)\right\}_{\textup{H}}: = \frac{1}{2} \left( H(j\omega) + H^*(j \omega)\right)$ denotes its Hermitian part.}} \ir{For a vector $x\in \mathbb{R}^n$ we denote by \hirot{$x^{(k)}$} its $k$-th element.}
\todoKL{\st{Notation $\text{diag}$ is not defined. Also, use the same notation $\text{diag}$  or ${diag}$ (with or without the text environment).} }

\mm{
Given a set $\mathcal{X}$, $\textrm{int}{\mathcal{X}}$ denotes the interior of $\mathcal{X}$ and $\textrm{bd}\mathcal{X}$ denotes the boundary of $\mathcal{X}$. The normal cone of $\mathcal{X}$ at a point $x\in \mathcal{X}$ is defined by
$N_{\mathcal{X}} (x) = \left\{ v: \langle v, y - x\rangle \leq 0, ~ \forall  y \in \mathcal{X} \right\}$. }

Given $0 \leq \mu \leq L$, we denote $S(\mu, L)$ as the set of functions $f: \mathbb{R}^{d} \rightarrow \mathbb{R}$ that are continuously differentiable, $\mu$-strongly convex and $L$-smooth, i.e., $\forall x , ~ y$,
 \begin{align*}
 	\mu \| x - y \|^2 \leq \left( \nabla f (x) - \nabla f(y) \right)^T(x - y) \leq L \| x - y \|^2.
 \end{align*}
In this work, we assume $\mu > 0$ for all the functions we study if not specified otherwise. The condition number $\kappa$ of functions in $S(\mu, L)$ is defined by $\kappa := L/\mu \geq 1$.

A memoryless nonlinearity $\psi : \mathbb{R} \to \mathbb{R}$, \revise{$\psi (0) = 0$,}
is said to be \textit{sloped-restricted} on sector $[\alpha,  \beta]$, if $\alpha \leq \frac{\psi (x) - \psi (y)}{x - y} \leq \beta $, for all $x \neq y$ and $x, y \in \mathbb{R}$.

\subsection{Integral quadratic constraints}
\ic{To facilitate readability some preliminaries on integral quadratic constraints (\cite{rantzer1996kalman}, \cite{jonsson2001lecture}) are included in this section.}

In continuous time, a bounded operator $\Delta : \mathbf{L}_2^{n} [0, \infty) \rightarrow \mathbf{L}_{2}^{m} [0, \infty)$ is said to satisfy the IQC defined by $\Pi$, denoted by $\Delta \in \text{IQC}(\Pi)$, if
\begin{align}\label{eq:IQC definition}
	\int_{-\infty}^{\infty} \begin{bmatrix} \hat{v} (j\omega) \\ \hat{w} (j \omega) \end{bmatrix}^*
	\Pi(j\omega)
	\begin{bmatrix} \hat{v} (j\omega) \\ \hat{w} (j \omega) \end{bmatrix} d\omega \geq 0
\end{align}
for all $v \in \mathbf{L}_{2}^{n} [0, \infty)$ and $w = \Delta(v)$, where $\hat{v}(j\omega)$, $\hat{w} (j\omega)$ are the Fourier transforms of $v$, $w$, respectively, and $\Pi (j\omega)$ can be any measurable Hermitian valued function.
In discrete time, condition \eqref{eq:IQC definition} reduces to
\begin{align*}
	\int_{- \pi}^{\pi} \begin{bmatrix} \hat{v} ( e^{j\omega} ) \\ \hat{w} (e^{j \omega}) \end{bmatrix}^*
	\Pi (e^{j \omega})
	\begin{bmatrix} \hat{v} (e^{j\omega} ) \\ \hat{w} ( e^{j \omega} ) \end{bmatrix}
	d \omega \geq 0
\end{align*}
for all $ v \in {l}_{2}^{n}(\mathbb{Z}_{+})$, and $w = \Delta (v)$, \lmm{where $\hat{v} ( e^{j\omega} )$, $\hat{w} ( e^{j\omega} )$ represent the discrete-time Fourier transform of $v$, $w$ respectively. }

Define the truncation operator $P_T$ which does not change a function on the interval $[0, T]$ and gives the value zero on $(T, \infty]$.
The operator $\Delta$ is said to be \textit{causal} if $P_T \Delta P_T = P_T \Delta$, for all $T \geq 0$.
Consider the interconnection
\begin{equation}\label{eq:feedback interconnection model}
\begin{aligned}
	v = & Gw + g\\
	w = & \Delta (v) + e
\end{aligned}
\end{equation}
where $g \in \mathbf{L}_{2e}^{l}[0, \infty)$, $e \in \mathbf{L}_{2e}^{m}[0, \infty)$, $G$ and $\Delta$ are two causal operators on $\mathbf{L}_{2e}^{m}[0, \infty)$, $\mathbf{L}_{2e}^{l}[0, \infty)$, respectively.
\lmm{It is assumed that $G$ is a stable linear time-invariant operator and $\Delta$ has bounded gain \lir{denoted as $\|\Delta\|$}.}
The feedback interconnection of $G$ and $\Delta$ is \textit{well-posed} if the map $(v, w) \mapsto (e, g)$ defined by \eqref{eq:feedback interconnection model} has a causal inverse on $\mathbf{L}_{2e}^{m+l}[0, \infty)$. The interconnection is \textit{stable} if, in addition, the inverse is bounded, i.e., there exists a constant $c > 0$ such that $\int_{0}^{T} \left( |v|^2  + |w|^2\right) d t \leq c \int_{0}^{T} \left( |g|^2 + |e|^2 \right) dt$.
System \eqref{eq:feedback interconnection model} with linear $G$ and \ml{static} nonlinear $\Delta$ is called the \textit{Lur'e system}.

\todoing{\st{I think Lure system is about static nonlinearities}}
We adopt the following IQC theorem for analysis.
\begin{theorem}[\hspace{1sp}\cite{megretski1997system}]\label{thm:IQC theorem}
	Let $G(s) \in \mathbf{RH}_{\infty}^{l \times m}$, and let $\Delta$ be a bounded causal operator. Assume that:
	\begin{enumerate}
		\item for every $\tau \in [0,1]$, the interconnection of $G$ and $\tau \Delta$ is well-posed;
		\item for every $\tau \in [0,1]$, the IQC defined by $\Pi$ is satisfied by $\tau \Delta$;
		\item there exists $\epsilon > 0$ such that
		\begin{align}\label{eq:IQC theorem condition 3}
			\begin{bmatrix}
				G (j \omega) \\ I
			\end{bmatrix}^{*}
			\Pi (j\omega)
			\begin{bmatrix}
				G (j \omega) \\ I
			\end{bmatrix}
			\leq - \epsilon I, ~ \forall \omega \in \mathbb{R}.
		\end{align}
	\end{enumerate}
	Then, the interconnection of $G$ and $\Delta$ is stable.
\end{theorem}

\begin{remark}
If $\Pi (j \omega) = \begin{bmatrix}
	\Pi_{11} (j \omega) & \Pi_{12}(j \omega) \\ \Pi_{12}^* (j \omega) & \Pi_{22} (j \omega)
\end{bmatrix}$ satisfies $\Pi_{11} (j \omega) \geq 0$ and $\Pi_{22} (j \omega) \leq 0$, then the condition $\Delta \in \text{IQC}(\Pi)$ implies that $\tau \Delta \in \text{IQC} (\Pi)$ for \ic{all $\tau \in [0, 1]$.} 
\hiro{
\lir{$\Pi (j\omega)$ in \cite{megretski1997system} is required to be essentially bounded, and extensions that allow to include also the Popov IQC were derived in \cite{jonsson1997stability}.}}
\end{remark}
The IQC theorem \il{for discrete-time} systems can be found in, e.g., \cite{jonsson2001lecture}.

\subsection{Mirror descent algorithm}
Consider the optimization problem
\begin{align}\label{eq:optimization problem}
\min_{x \in \mathcal{X}} f(x)
\end{align}
where $\mathcal{X}$ is \icl{a} closed and convex constraint set and $ \mathcal{X} \subseteq \mathbb{R}^{d}$, $f$ is the objective function and $f \in S(\mu, L)$.
We will consider the unconstrained case in Section~\ref{Continuous-time mirror descent method}, and \ref{Discrete-time mirror descent method} first, i.e., $\mathcal{X} = \mathbb{R}^{d}$,
and the constrained case in Section~\ref{Projected MD and convergence rate}.
%

\lmm{When $\mathcal{X}=\mathbb{R}^d$}, we can solve \eqref{eq:optimization problem} with the well-known gradient descent (GD) algorithm
$
x_{k+1} = x_k - \eta \nabla f(x_k),
$
or equivalently,
\begin{align*}
x_{k+1} = \underset{x \in \mathbb{R}^{d}}{\text{argmin}} \left\{ \nabla f(x_k)^T  x  + \frac{1}{2 \eta} \| x - x_{k} \|^2_2  \right\}
\end{align*}
where $\eta > 0$ is a fixed stepsize. Observe that the Euclidean norm used above can be replaced by other \icl{distance \il{measures}}
to generate new algorithms.

The \textit{Bregman divergence} defined with respect to a distance generating function (DGF) $\phi: \mathbb{R}^{d} \rightarrow \mathbb{R}$ is given by
\begin{align}\label{eq: Bregman divergence}
D_{\phi} (y, x) = \phi (y) - \phi (x) - ( y - x )^T \nabla \phi (x).
\end{align}
where $\phi(x) \in S (\mu_{\phi}, L_{\phi})$.
Then, the MD algorithm is given by
\begin{align}\label{eq:MD algorithm compact}
x_{k+1} = \underset{x \in \mathbb{R}^{d}}{\text{argmin}} \left\{ \nabla f(x_k)^T x + \frac{1}{\eta} D_{\phi} (x, x_{k}) \right\}.
\end{align}
\mm{When \revise{$\phi (x) = \frac{1}{2} \| x \|^2$},  the MD algorithm is reduced to the GD algorithm.}
Denote $\bar{\phi}$ as the convex conjugate of function $\phi$, i.e.,
\begin{align}\label{eq:convex conjugate}
\bar{\phi}(z) = \sup_{x \in \mathbb{R}^{d}} \left\{ x^T z - \phi (x) \right\}.
\end{align}
Denote $\mu_{\bar{\phi}} = \left(L_{\phi} \right)^{-1}$, and $L_{\bar{\phi}} =\left( \mu_{\phi}\right)^{-1}$. It follows
\revise{from \cite[Proposition 11.3 and 12.60]{rockafellar2009variational}}
that $\bar{\phi} \in S (\mu_{\bar{\phi}}, L_{\bar{\phi}} )$, and
$
z = \nabla \phi(x) \Longleftrightarrow x = \nabla \bar{\phi} (z).
$
In other words, $\nabla \bar{\phi}$ is the inverse function of $\nabla \phi$.
Then, from the optimality condition for the right-hand side of \eqref{eq:MD algorithm compact}, the MD algorithm can be written as
\begin{align*}
z_{k+1}  = z_{k} - \eta \nabla f (x_{k}), \quad x_{k+1} = \nabla \bar{\phi} (z_{k+1} )
\end{align*}
or equivalently,
\begin{align}\label{eq:MD discrete-time composition}
z_{k+1} = z_{k} - \eta ( \nabla f  \circ \nabla \bar{\phi} )(z_{k})
\end{align}
where $\circ$ represents \il{composition of functions}.
Similarly, the continuous-time MD algorithm can be given by
\begin{align}\label{eq:MD continuous-time composition}
\dot{z} \nkl{(t)} =  - \eta (\nabla f  \circ \nabla \bar{\phi})(z\nkl{(t)}).
\end{align}
\ml{Any equilibrium point of the above systems satisfies $\nabla f \left( \nabla \bar{\phi} (z^\text{\nkl{opt}}) \right) = \nabla f(x^\text{\nkl{opt}}) = 0_d$, which yields the optimal solution to problem \eqref{eq:optimization problem}.
}
\nkl{In the remainder of this paper, the time \icc{argument} 
in the continuous-time case
will be omitted to simplify the notation.}

\todoing{\st{how is $\phi$ chosen? Can it be arbitrary? Why is the equilibrium point unchanged?}}
\ml{Note that the DGF $\phi$ can be an arbitrary function in $S (\mu_{\phi}, L_{\phi})$.
\il{Function} $\phi$ is \il{usually} chosen such that its convex conjugate is easily computable.
\icr{One of the motivations} \icr{in the choice of $\phi$} is to generate a distance function that reflects the geometry of the given constraint set $\mathcal{X}$ so \iclr{that projections in $\mathcal{X}$ can be avoided\footnote{\iclr{In particular, projections in $\mathcal{X}$ are avoided when condition \eqref{eq:smooth assumption} stated in \ml{Section~}\ref{sec:ProjMD} is satisfied (\cite{beck2003mirror}).}}}
during \icr{the implementation of the algorithm}.
\revise{For example,  choosing the negative entropy $\phi (x) = \sum_{i = 1}^{d} x^{(i)} \ln x^{(i)} $ avoids the direct use of projection when considering minimization over the unit simplex $\mathcal{X} = \left\{x = [x^{(1)},\ldots, x^{(d)}]^T \in \mathbb{R}^{d}_{+}: \sum_{i = 1}^{d} x^{(i)} = 1\right\}$ (\cite{beck2003mirror,bubeck2015convex}).}
\revise{Additionally, the \iclr{choice of $\phi$ changes 
the geometry of gradient} \icr{descent}, potentially improving the convergence rate (\cite{maddison2021dual}).}

\section{Continuous-time mirror descent method}\label{Continuous-time mirror descent method}
\lmm{
\ic{We start 
our discussion in this paper} with the analysis of the continuous-time MD method. \ic{A main result in this section is to reveal a connection between the Bregman divergence and the Popov criterion when these are used as tools for stability analysis.  In particular, we show that the use of the Bregman divergence from the optimal solution as a Lyapunov function, is a special case of Lyapunov functions that follow from the Popov criterion when it is applied to an appropriate reformulation of the \revise{MD dynamics}. This connection is used as a basis to develop an IQC framework that can provide convergence rate bounds with reduced conservatism. In particular, this framework is exploited in subsequent sections to provide convergence rate bounds in a discrete-time setting and also when projections are considered to constrain the dynamics in prescribed sets. }}
%

\subsection{MD algorithm in the form of Lur'e systems}
\ml{The composition of operators in \eqref{eq:MD continuous-time composition} hinders the \il{direct} application of \il{an} IQC framework since the \icl{composite} operator may not belong to the original classes of the two operators, e.g., the composition of two monotone operators is not necessarily monotone.}
Nevertheless, the cascade connection of two nonlinear operators can be transformed into the feedback interconnection of a linear system with the direct sum of the two nonlinear operators, similarly to the example in \cite{megretski1997system}.
\mm{Thus, \eqref{eq:MD continuous-time composition} can be written into a Lur'e system
\begin{align*}
\dot{z}  =\hspace{-1mm} -\eta u_1', ~
 u_1' = \hspace{-0.5mm} \nabla f (y_1), ~
 y_1 = u_2',  ~
 u_2' =\hspace{-0.5mm}  \nabla \bar{\phi} (y_2), ~
 y_2 = z.
\end{align*}
\lmm{Moreover, we would like to transform the \ic{nonlinearities into sector-bounded ones that allow}
to generate IQCs that fulfil the \ic{second 
condition} in Theorem~\ref{thm:IQC theorem} (\ic{discussed in} Section~\ref{sbusection Sector IQC}).} Therefore,  the continuous-time MD algorithm \eqref{eq:MD continuous-time composition} is rewritten \ic{as} 
\begin{align}\label{eq:linear system continuous-time}
\dot{z}  = A z + B u, \quad
y =
C z +
D u
\end{align}
where $u = \begin{bmatrix} u_1\\ u_2 \end{bmatrix}$, $y = \begin{bmatrix}
y_1 \\ y_2
\end{bmatrix}$, the system matrices are
\begin{align}\label{eq:system matrices continuous-time}
\begin{bmatrix}
	\begin{array}{c|c}
      A & B\\
      \hline
      C & D
    \end{array}
\end{bmatrix}
=
\begin{bmatrix}
	\begin{array}{c|cc}
		-\eta \mu_f \mu_{\bar{\phi}} I_d &
		-\eta I_d & -\eta \mu_f I_d \\
		\hline
		\mu_{\bar{\phi}} I_d & 0_d & I_d\\
		I_d & 0_d & 0_d
		\end{array}
\end{bmatrix}
\end{align}
and the system input is
\begin{align}\label{eq:system input continuous-time}
\begin{bmatrix}
u_1 \\ u_2
\end{bmatrix}
= &
\begin{bmatrix}
\nabla f (y_1) - \mu_{f} y_1 \\ \nabla \bar{\phi} (y_2) - \mu_{\bar{\phi}} y_2
\end{bmatrix}.
\end{align}
The transfer function matrix of the linear system is
\begin{equation}\label{eq:transfer function continuous-time}
\begin{aligned}
G(s) = & C(sI_{d} - A)^{-1} B + D\\
= &
\frac{1}{s + \eta \mu_{f} \mu_{\bar{\phi}} }
\begin{bmatrix}
- \eta \mu_{\bar{\phi}} & s\\
- \eta & - \eta \mu_{f}
\end{bmatrix} \otimes I_d
\end{aligned}
\end{equation}
where $\otimes$ denotes the Kronecker product.
Next, \ic{we} define $z^\text{\nkl{opt}}$, $x^\text{\nkl{opt}}$ as the unique \nkl{optimal  state} with corresponding \ic{variables} $y^\text{\nkl{opt}}$ and $u^\text{\nkl{opt}}$.
Let $\tilde{z} = z - z^\text{\nkl{opt}}$, $\tilde{y} = y - y^\text{\nkl{opt}}$, $\tilde{u} = u - u^\text{\nkl{opt}}$.
We obtain the error system
\begin{equation}\label{eq:mirror descent algorithm continuous-time}
\dot{\tilde{z}} =   A \tilde{z} + B \tilde{u}, \quad
\tilde{y} =
C
\tilde{z} +
D
\tilde{u}\end{equation} with \begin{equation}
\begin{split}
\tilde{u}
:= & {\Delta} \left( \begin{bmatrix} y_1 - y^\text{\nkl{opt}}_1\\ y_2 - y^\text{\nkl{opt}}_2 \end{bmatrix} \right)
 =  \begin{bmatrix} \Delta_1 \left(y_1 - y^\text{\nkl{opt}}_1 \right) \\ \Delta_2 \left( y_2 - y^\text{\nkl{opt}}_2 \right)\end{bmatrix},
\end{split}
\label{eq:nonlinearity}
\end{equation}
where $\Delta_1 (x)$, $\Delta_2 (x)$ are defined by
$\Delta_1 (x)
= (\nabla f (x + y_1^{\textup{opt}})  -  \mu_{f} (x + y_1^{\textup{opt}}) ) - \left(\nabla f (y^\text{\nkl{opt}}_1) - \mu_{f} y^\text{\nkl{opt}}_1 \right)$,
$ \Delta_2 (x)
=  ( \nabla \bar{\phi} (x + y_2^{\textup{opt}})  - \mu_{\bar{\phi}} (x + y_2^{\textup{opt}}) ) - \left(\nabla \bar{\phi} (y^\text{\nkl{opt}}_2) - \mu_{\bar{\phi}} y^\text{\nkl{opt}}_2 \right).
$
\ml{The above error system is in the form of \il{a Lur'e system} \eqref{eq:feedback interconnection model}, where $v = \tilde{y}$, $w = \tilde{u}$, $e = 0_d$, and $g$ \il{is a trajectory that} represents the \il{effect of the} initial condition. \ic{The transformation is illustrated in} 
Fig.~\ref{fig:feedback interconnection of nonlinear functions}.}

\begin{figure}[htbp]
\centering
\subfigure[]{\includegraphics[width = 0.8\linewidth]{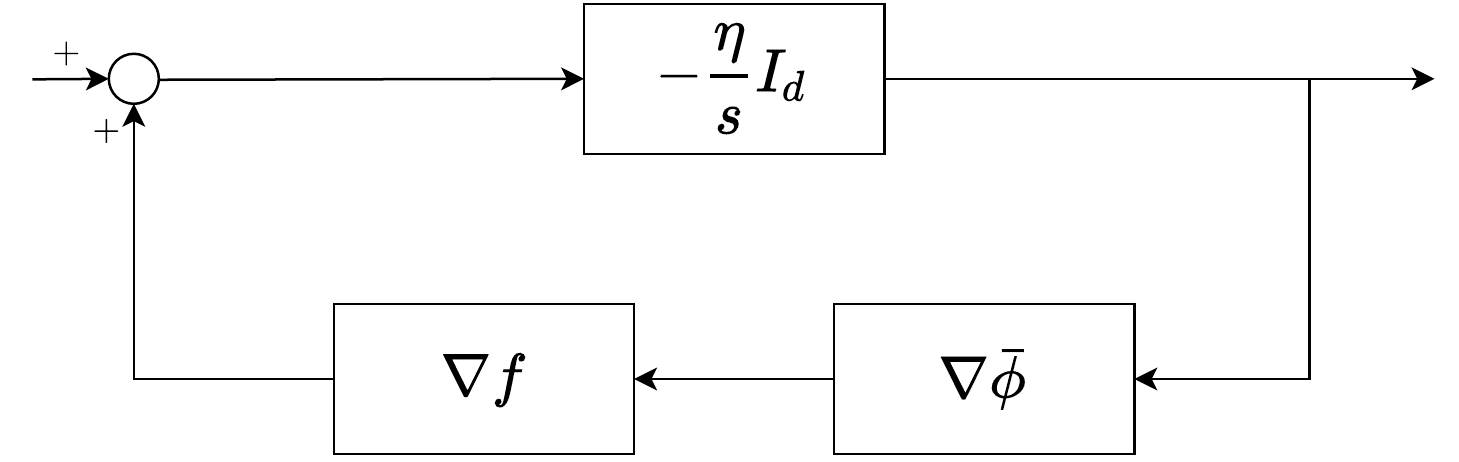}\label{fig: feedback cascade}}
\subfigure[]{\includegraphics[width = 0.8\linewidth]{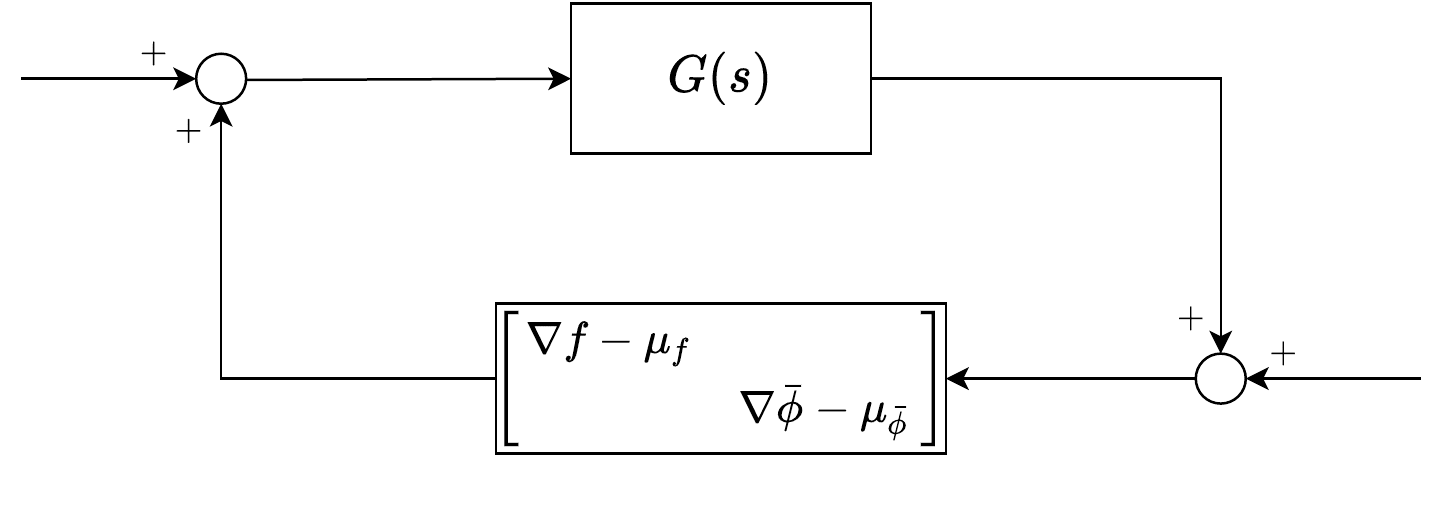}\label{fig: feedback_cascade_equivalence}}
\subfigure[]{\includegraphics[width = 0.8\linewidth]{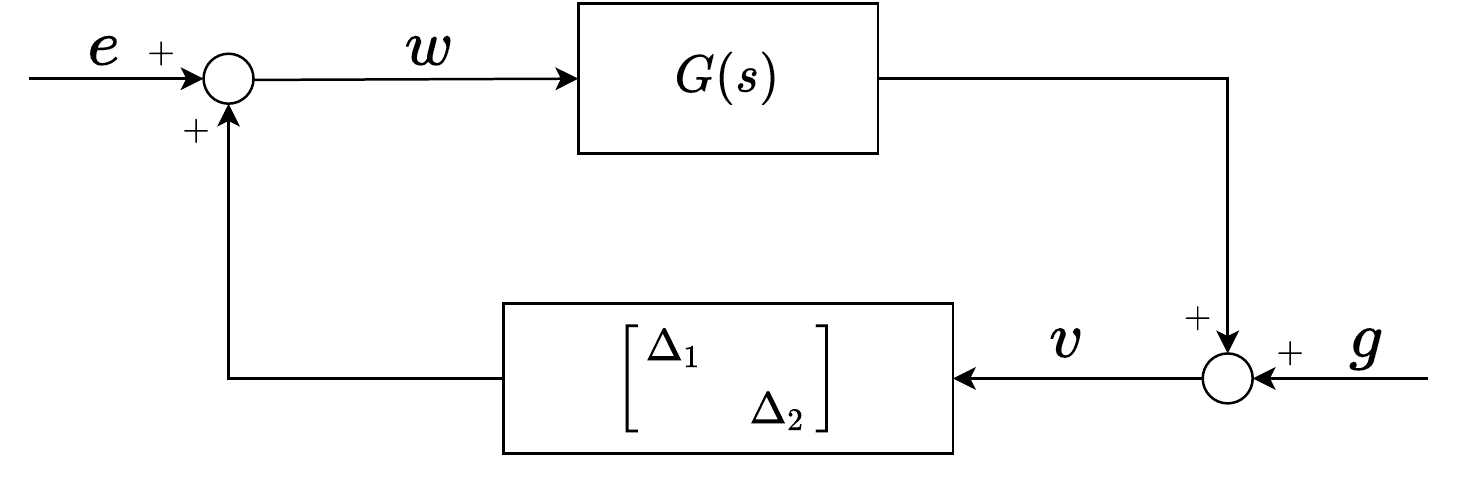}\label{fig: error_system}}
\caption{\ml{Transformation of the MD method to \il{a Lur'e} system. (a) represents the composition of operators, which is transformed to the direct sum of operators in (b), where $G(s)$ is given by \eqref{eq:transfer function continuous-time}.
(c) is the error system in \eqref{eq:mirror descent algorithm continuous-time}, and $\Delta_1$, $\Delta_2$ are given by \eqref{eq:nonlinearity}. }}
\label{fig:feedback interconnection of nonlinear functions}
\end{figure}

\subsection{IQCs for gradients of convex functions}
\todoiny{\st{I think more text is needed before starting section III.B to justify why we are bringing all these materials.}}
\ml{In this subsection, we will include a \icc{number} 
of useful IQCs for \icc{the} gradients of convex functions \icc{so as to} characterize the nonlinearity $\Delta$.
Note that conic combinations of various IQCs are also valid IQCs which better characterize the nonlinearity and lead to less conservative stability margins.
\todoiny{\st{Say something about the conic combination}}
}
\subsubsection{Sector IQC}\label{sbusection Sector IQC}
\ml{The Sector IQC \ic{stated in the following lemma is a} 
result of the co-coercivity of gradients.}
\begin{lemma}[\hspace{1sp}\cite{lessard2016analysis}]\label{lem:co-coercivity sector IQC}
Suppose $f \in S(\mu_f, L_f)$. For all $x, y$, the following quadratic constraint (QC) is satisfied,
\begin{align*}
\left[
\begin{smallmatrix}
y - x \\
\nabla f (y) - \nabla f(x)
\end{smallmatrix}
\right]^T
\left[
\begin{smallmatrix}
-2 \mu_f L_f I_d & (L_f + \mu_f) I_d \\
(L_f + \mu_f )I_d & -2 I_d
\end{smallmatrix}
\right]
\left[
\begin{smallmatrix}
y - x \\
\nabla f (y) - \nabla f(x)
\end{smallmatrix}
\right]
\geq 0.
\end{align*}
\end{lemma}

Note that as  $f \in S(\mu_f, L_f)$, $\bar{\phi} \in S(\mu_{\bar{\phi}}, L_{\bar{\phi}})$, \ml{then $f (\cdot) - \frac{1}{2} \mu_f \| \cdot \|^2 \in S ( 0, L_f - \mu_f )$ and $\bar{\phi} (\cdot) - \frac{1}{2} \mu_{\bar{\phi}} \|\cdot \|^2 \in S ( 0, L_{\bar{\phi}}- \mu_{\bar{\phi}} )$.}
\todoinc{are you meant to substract just a constant to the functions? what is the role of $\|y_i\|^2$, and the division by 2?}
Moreover, using Lemma~\ref{lem:co-coercivity sector IQC},  we have that,  given $\Delta$ defined in \ml{\eqref{eq:nonlinearity}}, $\Delta \in \text{IQC}( \Pi_{s})$, where
\begin{align}\label{eq:IQC sector bounded}
\Pi_{s} \hspace{-1mm} =
\left[
\begin{smallmatrix}
0_d & 0_d & \alpha_1 ( L_f - \mu_f )  I_d & 0_d\\
0_d & 0_d & 0_d & \alpha_2 ( L_{\bar{\phi}} - \mu_{\bar{\phi}} ) I_d\\
\alpha_1 ( L_f - \mu_f )  I_d & 0_d & -2 \alpha_1 I_d & 0_d\\
0_d & \alpha_2 ( L_{\bar{\phi}} - \mu_{\bar{\phi}} ) I_d & 0_d & -2 \alpha_2 I_d
\end{smallmatrix}
\right]
\end{align}
for some $\alpha_1, \alpha_2 \geq 0$.
\lmm{The second \ic{condition in Theorem~\ref{thm:IQC theorem} is always} 
satisfied with this IQC.}

\subsubsection{Popov IQC}\label{sec:PopovIQC}
\ml{The Popov IQC \ic{can be stated} 
as follows.}
\begin{lemma}[\hspace{1sp}\cite{jonsson1997stability}]\label{lem:popov IQC}
Suppose $f \in S(0, L)$. The nonlinearity $\nabla f(x) - \nabla f(x^\textup{\nkl{opt}})$ satisfies the Popov IQC 
\irr{defined \ir{by the}
multiplier}
\begin{align*}
\Pi_{P} (j \omega) =
	\pm \begin{bmatrix}
			0_d & - j \omega I _d \\ j\omega I_d & 0_d
		\end{bmatrix}
\end{align*}
\irr{in the sense described in \cite[Definition 1]{jonsson1997stability}}.
\end{lemma}
As $f (\cdot) - \frac{1}{2} \mu_f \| \cdot \|^2 \in S ( 0, L_f - \mu_f )$ and $\bar{\phi} (\cdot) - \frac{1}{2} \mu_{\bar{\phi}} \| \cdot \|^2 \in S ( 0, L_{\bar{\phi}}- \mu_{\bar{\phi}} )$,
using Lemma~\ref{lem:popov IQC}, we have $\Delta \in \text{IQC} (\Pi_{p})$, where $\Delta$ is defined in \eqref{eq:nonlinearity} and
\begin{align}\label{eq:Popov IQC for delta}
\Pi_{p}(j \omega) =
\left[
\begin{smallmatrix}
0_d & 0_d & -j \omega \beta_1  I_d & 0_d\\
0_d & 0_d & 0_d & - j \omega \beta_2 I_d\\
 j \omega \beta_1  I_d & 0_d & 0_d & 0_d\\
0_d & j \omega \beta_2 I_d & 0_d & 0_d
\end{smallmatrix}
\right]
\end{align}
for some $\beta_1, \beta_2 \in \mathbb{R}$.

\hiro{It is worth noting that the nonlinearity $\nabla f(x) - \nabla f(x^\textup{opt} )$ with $f \in S( 0, L)$ satisfies the IQC defined by \kl{$\Pi (j \omega)$ given by}}
\begin{align}\label{eq:Zames-Falb multipliers}
		\Pi (j \omega) = \begin{bmatrix}
			0_d & M^*(j\omega) I_d\\ M (j \omega) I_d & 0_d
		\end{bmatrix}
\end{align}
where $M(j\omega) = 1 - H(j\omega)$ and $H(j\omega)$ is the Fourier transform of any signal $h(t)$ such that $h(t) \geq 0$ and $\int_{-\infty}^{+\infty} |h(t)| dt \leq 1$.
\hiro{These $M(j \omega)$ are known as the \textit{Zames-Falb-O'Shea multipliers}.}
{\icl{It is well-known that the Zames-Falb-O'Shea multipliers can be used to \il{formulate} a wide class of IQCs that are satisfied by slope-restricted nonlinearities} \icl{(see also the discrete-time analogues in \cite{lessard2016analysis} formulated in time domain).}}

\lmm{
\begin{remark}
The Popov IQC can be obtained by letting $M(j\omega) = j\omega$ in \eqref{eq:Zames-Falb multipliers}. Note that \kl{the Popov multipliers} are not Zames-Falb-O'shea multipliers since the resulting $H (j \omega) = 1 - M(j \omega)$ is unbounded. Nevertheless, \icl{the} Popov multipliers can be treated as a \ml{limit} to the first order Zames-Falb-O'shea multipliers\footnote{\icr{We will show later on that conic combinations of Sector and Popov IQCs can recover the Bregman divergence based Lyapunov function in continuous time. The discrete-time \iclr{analogues} of Zames-Falb-O'Shea multipliers could be
beneficial in discrete time to reduce the conservatism, or when projections are additionally included.}} (\cite{carrasco2013equivalence}).
\end{remark}
}
\todoing{\st{without the boundedness assumption on the multiplier?}
Answer: As we discussed previously, the Popov $\Pi_{P}$ does not have boundedness assumption.}

\todoing{\st{if you have multiple IQCs their conic combination is also an IQC.}

Answer: This is now mentioned at the beginning of the subsection.}

\subsection{Convergence analysis via \il{IQCs} in \nkl{frequency domain}}
\todoing{\st{Need some introductory text to summarize the convergence results that will be presented.}}

In this subsection, we will present \ic{a 
convergence analysis for} 
the continuous-time MD method.
There exists rich literature showing the convergence of the MD method, e.g., \cite{nemirovskij1983problem,beck2003mirror,krichene2015accelerated}.
\lmm{Our approach using an IQC analysis framework \ic{recovers known 
conclusions} of convergence, but also lays the foundation for the subsequent \ic{sections where tight bounds on the convergence rate are obtained in discrete time. Furthermore, it allows to consider the more involved problem where projections to a convex set are additionally included.}} 

\begin{theorem}\label{thm:convergence continuous-time IQC}
The Lur'e system \icl{described by} \eqref{eq:feedback interconnection model} where $g \in \mathbf{L}_{2}[0, \infty)$, $e \in \mathbf{L}_{2}[0, \infty)$, $G(s)$ is given by \eqref{eq:transfer function continuous-time}, \ml{$\Delta$ is defined in \eqref{eq:nonlinearity}} with $f \in S(\mu_f, L_f)$, $\phi \in S(\mu_\phi, L_\phi)$, is stable.
\lmm{Furthermore, the trajectory of $x = \nabla \bar{\phi} (z)$ with any 
\il{initial} condition $z(0) = z_0$ of the MD algorithm \eqref{eq:MD continuous-time composition} converges to the optimal solution of problem \eqref{eq:optimization problem}.}
\todoing{\st{Need to be careful about the theorem statement - stability in IQCs is in an input/output sense, whereas here you want to deduce convergence for any initial conditions; see e.g. how Lessard paper handles this.}}
\end{theorem}
\revise{\begin{pf}
See Appendix~\ref{Proof of thm:convergence continuous-time IQC}.
\end{pf}}
\ml{\irr{The proof of Theorem~\ref{thm:convergence continuous-time IQC} can} be seen as an application of the multivariable Popov criterion,
which is \il{adapted below from} 
 \cite{moore1968generalization,khalil1996nonlinear,jonsson1997stability,carrasco2013equivalence}.
\begin{lemma}\label{lem:popov}
\nkl{Let $H (s) \in \mathbf{RH}_{\infty}^{p\times p}$ and let $\psi  : \mathbb{R}^{p} \to \mathbb{R}^{p}$, \revise{$\psi (0) = 0$},  be a memoryless nonlinearity composed of $p$ memoryless nonlinearities $\psi_i$ with each being }\ml{slope-restricted on sector \textup{[0,~$k_i$]}}, i.e., \ml{$0 \leq \frac{\psi_i(x_1) - \psi_i(x_2)}{x_1 - x_2} \leq k_i$, $\forall x_1 \neq x_2$}, $0< k_i \leq \infty$, for $1\leq i \leq p$.
If there exist constants $q_i \geq 0$ and $\gamma_i \geq 0$ such that $\left\{ Q K^{-1} + (Q + j \omega \Gamma) H (j\omega ) \right\}_{\hiro{\textup{H}}} \geq \delta\revise{I}$, $\forall \omega \in \mathbb{R}$, for some $\delta > 0$, where $K = \textup{diag}(k_1, \ldots, k_p)$, $Q = \textup{diag} (q_1, \ldots, q_p)$, $\Gamma = \textup{diag}(\gamma_1, \ldots, \gamma_p)$, and $\Gamma H(\infty) = 0$, then the negative feedback interconnection of $H(s)$ and $\psi$ is stable.
\end{lemma}
\begin{remark}
The third condition in Theorem~\ref{thm:IQC theorem}, \il{with the IQC used in the proof of Theorem~\ref{thm:convergence continuous-time IQC},}  is equivalent to the \lir{frequency domain} inequality condition \il{in Lemma~\ref{lem:popov} \hiro{with $H(s) = -G(s)$.}}
The consideration of different parameters $q_i$ is crucial since it provides \kl{more flexibility and thus less \ic{conservative}} results for the MIMO \il{case (\cite{moore1968generalization}).}
\hiro{\revise{The original Popov criterion requires that the stable linear system $H(s)$ is strictly proper (\cite{moore1968generalization,khalil1996nonlinear}). This restriction is relaxed to $\Gamma H(\infty) = 0$ in \cite{jonsson1997stability}.
}}
\todoing{\st{as above $\mathbf{L}_2$ stable needs clarification}}
\end{remark}

Theorem~\ref{thm:convergence continuous-time IQC} is based on conditions in frequency domain, which \icl{do not} describe the convergence rate of the MD algorithm. To this end, we will \icl{investigate, in the next subsection,} the MD method in time domain.
\subsection{\nkl{Convergence \hiro{rate} analysis via \il{IQCs} in time domain}}

In this subsection, we analyze convergence rate of the MD algorithm. Moreover,
we show that the Bregman divergence function, \icl{which is widely used as a Lyapunov function for the MD algorithm, is a special case of Lyapunov functions that are associated with the Popov criterion.} \todoing{this sentence is vague}

%
\todoiny{\st{More details are needed.}}

We can combine Theorem~\ref{thm:convergence continuous-time IQC} with what is proposed in \cite{hu2016exponential} \iclr{to} 
obtain a condition \revise{\iclr{that leads to lower} bounds on the} exponential convergence rate for the continuous-time MD method.
%
\revise{A continuous-time signal $x(t)$ converges to $x^{\textup{opt}}$ exponentially with rate $\rho \hiro{>} 0$ if there exists $c>0$ such that $\|x(t) - x^{\textup{opt}} \| \leq c e^{-\rho t} \| x(0) - x^{\textup{opt}}\|$, $\forall t \geq 0$.}
\begin{theorem}\label{thm:convergence rate continuous-time}
The continuous-time MD algorithm \eqref{eq:MD continuous-time composition} with $f \in S(\mu_f, L_f)$, $\phi \in S(\mu_\phi, L_\phi)$, converges exponentially to the optimal solution \nkl{with \ml{rate} $\rho$}
if there exist $P = P^T > 0$, $\Gamma = \textup{diag}(0, \gamma) \geq 0$, and $Q = \textup{diag} ( q_1, q_2) \geq 0$ such that
\begin{align}\label{eq:LMI exponential rate continuous-time}
	\begin{bmatrix}
		P\tilde{A}_{\hiro{\rho}} + \tilde{A}_{\hiro{\rho}}^TP & P \tilde{B} - \tilde{C}_{\hiro{\rho}}^T \\
		* & - \left( \tilde{D} + \tilde{D}^T \right)
	\end{bmatrix}
	\leq 0
\end{align}
where $\tilde{A}_{\hiro{\rho}} = A \hiro{+ \rho I}$, $\tilde{B} = - B$, $\tilde{C}_{\hiro{\rho}} = \mm{(Q + \hiro{2} \rho \Gamma)} C + \Gamma C A$, $\tilde{D} = - Q D + Q K^{-1} - \Gamma C B$, \hiro{$K = \textup{diag} (L_f-\mu_{f}, L_{\bar{\phi}} - \mu_{\bar{\phi}})$}, and $(A,B,C,D)$ is defined in \eqref{eq:system matrices continuous-time}.
Moreover, the
\irr{largest} lower bound on the convergence rate for any $f \in S(\mu_f, L_f)$, $\phi \in S(\mu_\phi, L_\phi)$ is given by $\rho = \eta \mu_f \mu_{\bar{\phi}}$.
\end{theorem}
\revise{\begin{pf}
See Appendix~\ref{Proof of thm:convergence rate continuous-time}.
\end{pf}}
\revise{\begin{remark}
\hu{It can be \hiro{easily} verified that \lrrr{there exist quadratic functions $f$, $\phi$ in the classes specified in the Theorem for which the rate bound $\rho = \eta \mu_f \mu_{\bar{\phi}}$ is tight.}}
\end{remark}}

\hiro{Consider the class of Lyapunov function \irr{candidates} \eqref{eq:lyapunov function popov 0} 
\irr{with $\rho=0$},
that is,}
\hiro{\begin{align}\label{eq:lyapunov function popov}
V = \frac{1}{2} \tilde{z}^T P \tilde{z} + \gamma \int_{0}^{y_2 - y^\textup{opt}_2} \psi_2 (\tau)  d \tau.
\end{align}
where $\psi_2 (\tau) = \Delta_2(\tau)$ \lr{with $\Delta_2$} \hu{defined} in \eqref{eq:nonlinearity}, \irr{and $P$ satisfies LMI \eqref{eq:LMI exponential rate continuous-time} with a strict inequality}.}
\irr{Such Lyapunov functions follow from the Popov criterion, i.e., it is known that when the conditions in the Popov criterion (Lemma \ref{lem:popov}) hold such a $P$ exists and \eqref{eq:lyapunov function popov} is a valid Lyapunov function for the  Lur'e system \eqref{eq:mirror descent algorithm continuous-time} \eqref{eq:nonlinearity} (\cite{khalil1996nonlinear}).}
\iclr{The following proposition establishes a link between the Bregman divergence and
the Popov criterion  in the context of the MD method.}
\revise{
\icr{\begin{proposition} {\bf(link between Popov criterion and Bregman divergence)} \label{Prop:Berg_Popov}
Consider the class of Lyapunov functions in \hiro{\eqref{eq:lyapunov function popov} \irr{for the MD method \eqref{eq:MD continuous-time composition}},
which follow from the Popov criterion.}
These Lyapunov functions include the Bregman divergence function $D_{\bar{\phi}} (z, z^{\textup{opt}})$ as a special case.  
\end{proposition}}}
\begin{pf}
As the nonlinear function \icr{$\psi_2$} is given by $\psi_2 (\tilde{y}_2) = \left( \nabla \bar{\phi} (\tilde{y}_2 + y_2^{\textup{opt}}) - \mu_{\bar{\phi}} (\tilde{y}_2 + y_2^{\textup{opt}}) \right) - \left(  \nabla \bar{\phi} (y^\textup{\nkl{opt}}_2) - \mu_{\bar{\phi}} y^\textup{\nkl{opt}}_2 \right)$,  \eqref{eq:lyapunov function popov} becomes
\begin{align}\label{eq:lyapunov function popov to bregman}
	V \hspace{-1mm} = & \frac{1}{2} \tilde{z}^T P \tilde{z} \hspace{-0.5mm} + \hspace{-0.5mm} \gamma \hspace{-1mm} \int_{y^\text{opt}_2}^{y_2} \hspace{-1mm} \left( \nabla \bar{\phi} (\tau) \hspace{-1mm} - \hspace{-1mm}\mu_{\bar{\phi}} \tau \right) \hspace{-1mm} - \hspace{-1mm} \left( \nabla \bar{\phi} (y^\text{\nkl{opt}}_2) \hspace{-1mm} - \hspace{-1mm} \mu_{\bar{\phi}} y^\text{\nkl{opt}}_2 \right) d \tau \nonumber\\
	= & \frac{1}{2} \tilde{z}^T P \tilde{z} + \gamma \hspace{-1mm} \int_{z^\text{\nkl{opt}}}^{z} \hspace{-1mm} \left( \nabla \bar{\phi} (\tau) \hspace{-1mm} - \hspace{-1mm}\mu_{\bar{\phi}} \tau \right) \hspace{-1mm} - \hspace{-1mm} \left( \nabla \bar{\phi} (z^\text{\nkl{opt}}) \hspace{-1mm} - \hspace{-1mm} \mu_{\bar{\phi}} z^\text{\nkl{opt}} \right) d \tau \nonumber\\
	= & \frac{1}{2} \tilde{z}^T P\tilde{z} + \gamma D_{\bar{\phi}}(z, z^\text{\nkl{opt}}) - \frac{\gamma \mu_{\bar{\phi}}}{2} \| \tilde{z} \|^2.
\end{align}
\hu{
\lr{We note that the LMI \eqref{eq: LMI with solution} in Appendix~\ref{Proof of thm:convergence rate continuous-time}, which is  an expanded version of \eqref{eq:LMI exponential rate continuous-time} when $P=pI_d>0$, $p\in\mathbb{R}$, \lr{holds with a strict inequality when}} $\rho = 0$, $P = \gamma \mu_{\bar{\phi}} I_d$, $q_{1} = \eta \gamma$, $q_{2} = 2 \eta \gamma \mu_{f} \mu_{\bar{\phi}}$, and $\eta > 0$. \lr{This} ensures that $V$ with this $P$ is a valid Lyapunov function, \lr{as follows from the proof of Theorem~\ref{thm:convergence rate continuous-time}.}
}
In addition, let $\gamma = 1$,  then the Lyapunov function \eqref{eq:lyapunov function popov} for the Popov criterion reduces to the Bregman divergence from the optimal solution,  $D_{\bar{\phi}}(z, z^\text{opt})$.
\end{pf}
\begin{remark}
\icr{Proposition \ref{Prop:Berg_Popov} shows that
the Bregman divergence from the optimal solution, which is the Lyapunov function commonly used in the literature for the MD algorithm (\cite{nemirovskij1983problem, krichene2015accelerated}), is a special case of Lyapunov functions that follow from the Popov criterion.
This establishes an important connection between the Bregman divergence and the Popov criterion, which illustrates the significance of an IQC analysis in the reformulated problem that has been considered.
In particular, an IQC analysis allows multiple IQCs to be combined via conic combinations, and can therefore potentially lead to less conservative convergence rate bounds.
This becomes relevant in the discrete-time case and the case where projections are additionally included, which are addressed subsequently in this work.}
\end{remark}

\section{Discrete-time mirror descent method}\label{Discrete-time mirror descent method}
\lmm{Having observed the effectiveness of IQCs in continuous time, we now turn our attention to the convergence rate analysis of MD method \ic{in discrete} time.}
\subsection{MD algorithm in the form of Lur'e systems}
Similar to the continuous-time case, the discrete-time MD \il{algorithm} in \eqref{eq:MD discrete-time composition} can be rewritten \ic{as} 
the following Lur'e system,
\begin{align}\label{eq:mirror descent discrete-time}
	z_{k+1} = A z_{k} + B u_k, \quad y_k = C z_k + D u_k
\end{align}
where $ u_k = \begin{bmatrix} u_k^{(1)} \\ u_k^{(2)} \end{bmatrix}$, $y_k = \begin{bmatrix} y_k^{(1)}\\ y_k^{(2)} \end{bmatrix}$ the system matrices are
\begin{align}\label{eq:system matrices discrete-time}
\begin{bmatrix}
	\begin{array}{c|c}
      A & B\\
      \hline
      C & D
    \end{array}
\end{bmatrix}
=
\begin{bmatrix}
	\begin{array}{c|cc}
		(1-\eta \mu_f \mu_{\bar{\phi}})I_d &
		-\eta I_d & -\eta \mu_f I_d \\
		\hline
		\mu_{\bar{\phi}} I_d & 0_d & I_d\\
		I_d & 0_d & 0_d
		\end{array}
\end{bmatrix}
\end{align}
and the system input is
\begin{align}\label{eq:input discrete-time}
	\begin{bmatrix}
u_k^{(1)} \\ u_k^{(2)}
\end{bmatrix}
= &
\begin{bmatrix}
\nabla f (y_k^{(1)}) - \mu_{f} y_k^{(1)} \\ \nabla \bar{\phi} (y_k^{(2)}) - \mu_{\bar{\phi}} y_k^{(2)}
\end{bmatrix}
.
\end{align}
\lmm{We denote} \ml{$z^\textup{\nkl{opt}}$ as the optimal \il{value of $z$ at} steady state, with corresponding \il{equilibrium values} $y^\textup{\nkl{opt}} = \left( y^{(1),\nkl{\text{opt}}}, y^{(2),\nkl{\textup{opt}}}\right)$, $x^\textup{\nkl{opt}}$, and $u^\textup{\nkl{opt}}$. \ic{Defining} $\tilde{u}_{k} = u_k - u^\textup{\nkl{opt}}$, \ic{we have}}
\begin{align}\label{eq:nonlinearity discrete-time}
\tilde{u}_{k}
= \Delta
\left(
\begin{bmatrix}
	y_k^{(1)} - y^{(1),\nkl{\textup{opt}}} \\ y_k^{(2)} - y^{(2),\nkl{\textup{opt}}}
\end{bmatrix}
\right)
\end{align}
\ml{where the nonlinear operator $\Delta$ in \eqref{eq:nonlinearity discrete-time} is \icl{the same as that used for the continuous-time algorithm in \eqref{eq:nonlinearity}.}}
%

\subsection{Convergence rate \ic{bounds via IQCs}}
\hiro{There is no exact counterpart for the Popov criterion in \hiro{discrete time,}
though we could derive LMI conditions for the discrete-time system \eqref{eq:system matrices discrete-time}, \eqref{eq:nonlinearity discrete-time} via \irrr{related} Lyapunov functions (\cite{park2019less}).
Similar ones are the discrete-time Jury-Lee criteria (\cite{jury1964stability,haddad1994parameter}).}
\ml{\hiro{Nonetheless,} we remark that in discrete time, all IQCs to characterize monotone and bounded nonlinearities are within the set of Zames-Falb-O'Shea IQCs \li{(\cite{carrasco2016zames})}.}
\todoing{\st{do not fully understand the previous sentence. What do you mean by saying the multipliers "preserve the positivity of monotone and bounded nonlinearities".}}
Therefore, we can directly apply the class of 
Zames-Falb-O'Shea \il{IQCs with a state-space representation as in} \cite{lessard2016analysis}. We will only adopt a simple type of the Zames-Falb-O'Shea IQC here, \ml{as \il{this} is sufficient to obtain a tight convergence rate for the unconstrained MD method.}
\todoing{\st{is this done for simplicity in the presentation?}\\
Same as the last answer. The simplest type is sufficient to obtain a tight convergence rate.
}

\begin{lemma}[\hspace{1sp}\cite{lessard2016analysis}]\label{lem:weighted-off-by-one IQC}
	Suppose $f \in S(\mu_f, L_f)$, $x^\textup{opt}$ is a reference point. The nonlinearity $\nabla f(x) - \nabla f(x^\textup{opt} )$ satisfies the \lmm{weighted off-by-one} IQC defined by $\Pi_{w, f} (j \omega)$ given by
	\begin{align}
		\Pi_{w, f} \hspace{-1mm} = & \Psi_{w,f}^{*} M \Psi_{w,f}, \quad M \hspace{-1mm} = \hspace{-1mm} \begin{bmatrix} 0_{d} & I_{d} \\ I_{d} & 0_{d}\end{bmatrix},\end{align}
\ml{where $\Psi_{w,f}$ is a transfer function matrix with the following state-space representation,}
	\begin{align}\label{eq:weighted-off-by-one iqc state space}
	\hspace{-2mm}
	\Psi_{w, f}:
	\begin{bmatrix}
		\hspace{-0.5mm}
			\begin{array}{c|cc}
				A_{\Psi_{w,f}} & B_{\Psi_{w,f}}^{y} & B_{\Psi_{w,f}}^{u} \\
				\hline
				C_{\Psi_{w,f}} & D_{\Psi_{w,f}}^{y} & D_{\Psi_{w,f}}^{u}
			\end{array}
		\hspace{-0.5mm}
		\end{bmatrix} \hspace{-0.5mm} \ic{:=} \hspace{-0.5mm}
		\begin{bmatrix}
		\hspace{-0.5mm}
			\begin{array}{c|cc}
				0_{d} & - K I_{d} & I_{d} \\
				\hline
				\bar{\rho}^2 I_{d} & K I_{d} & -I_{d}\\
				0_{d} & 0_{d} & I_{d}
			\end{array}
		\hspace{-0.5mm}
		\end{bmatrix} \hspace{-1mm}
	\end{align}
\ic{where} $K = L_f - \mu_f \geq 0$ and $\bar{\rho} \geq 0$.
\end{lemma}
\lmm{For $\bar{\phi} \in S ( \mu_{\bar{\phi}}, L_{\bar{\phi}}) $,  we define $\Psi_{w, \bar{\phi}}$ similarly to \eqref{eq:weighted-off-by-one iqc state space} with matrices $A_{\Psi_{w,\bar{\phi}}}, \ldots, D_{\Psi_{w,\bar{\phi}}}^{u}$.}

From Lemma~\ref{lem:co-coercivity sector IQC}, we can obtain that $\nabla f (x) - \nabla f(x^{\textup{opt}})$ satisfies the IQC defined by
\begin{align}
	& \Pi_{s, f} = \Psi_{s,f}^* M \Psi_{s,f}, \quad M = \begin{bmatrix} 0_{d} & I_{d} \\ I_{d} & 0_{d}\end{bmatrix}
\end{align}
\ml{where $\Psi_{s}$ is a transfer function matrix with the following state-space representation,}
\begin{align}\label{eq:sector iqc state space}
\Psi_{s, f}:
\begin{bmatrix}
	\begin{array}{c|cc}
		A_{\Psi_{s,f}} & B_{\Psi_{s,f}}^{y} & B_{\Psi_{s,f}}^{u}\\
		\hline
		C_{\Psi_{s,f}} & D_{\Psi_{s}}^{y} & D_{\Psi_{s,f}}^{u}
		\end{array}
\end{bmatrix} \ic{:=}
	\begin{bmatrix}
	\begin{array}{c|cc}
		0_{d} & 0_{d} & 0_{d} \\
		\hline
		0_{d} & K I_{d} & -I_{d}\\
		0_{d} & 0_{d} & I_{d}
		\end{array}
\end{bmatrix}.
\end{align}
\ic{where} $K = L_f - \mu_f \geq 0$. \lmm{For $\bar{\phi} \in S ( \mu_{\bar{\phi}}, L_{\bar{\phi}}) $,  we define $\Psi_{s, \bar{\phi}}$ similarly to \eqref{eq:sector iqc state space} with matrices $A_{\Psi_{s,\bar{\phi}}}, \ldots, D_{\Psi_{s,\bar{\phi}}}^{u}$.}

Then, we can characterize \revise{upper bounds on the} convergence rate for the discrete-time MD algorithm by applying the discrete-time IQC theorem \revise{with the IQC constructed from the conic combination of a sector
and \mmli{a weighted} off-by-one IQC for both nonlinearities in \eqref{eq:nonlinearity discrete-time}}.
\revise{A sequence $\{x(k)\}$ converges to $x^{\textup{opt}}$ linearly with rate $\rho \in (0, 1)$ if there exists $c >0$ such that $\| x(k) - x^{\textup{opt}}\| \leq c \rho^{k} \|x (0) - x^{\textup{opt}}\|$, for all $k \in \mathbb{N}$.}
\begin{theorem}\label{thm:convergence rate discrete-time}
	The discrete-time MD algorithm \eqref{eq:MD discrete-time composition} with $f \in S(\mu_f, L_f)$ and $\phi \in S(\mu_\phi, L_\phi)$ converges linearly \il{with} rate $\bar{\rho} \leq \rho \leq 1$ if the following LMI is feasible for some $P = P^T > 0$, $\alpha = \textup{diag} \{\alpha_1, \ldots, \alpha_4\} \geq 0$,
\begin{align}\label{eq:LMI IQC discrete-time}
	\left[
	\begin{smallmatrix}
		\hat{A}^T P \hat{A} - \rho^2 P & \hat{A}^T P \hat{B}\\
		* & \hat{B}^T P \hat{B}
	\end{smallmatrix} \right]
	+
	\left[\begin{smallmatrix}
		\hat{C} & \hat{D}
	\end{smallmatrix} \right]^T
	\left( \alpha \otimes M \right)
	\left[
	\begin{smallmatrix}
		\hat{C} & \hat{D}
	\end{smallmatrix} \right]
	\leq 0
	\end{align}
where
\hiro{\begin{align}
	\hat{A} = &
	\begin{bmatrix}
		A & 0_{d \times d}  &  0_{d \times d}\\
		B_{{\Psi_{w,f}}}^{y} C_1  &  A_{\Psi_{w,f}} & 0_{d \times d}\\
	 	B_{{\Psi_{w, \bar{\phi}}}}^{y} C_2 & 0_{d \times d} & A_{\Psi_{w,\bar{\phi}}}
	\end{bmatrix}, \nonumber \\
	\hat{B} = &
	\begin{bmatrix}
    B_1 & B_2\\
    B_{\Psi_{w, f}}^{y} D_{11} + B_{\Psi_{w, f}}^{u} & B_{\Psi_{w, f}}^{y} D_{12}\\
    B_{\Psi_{w, \bar{\phi}}}^{y} D_{21} & B_{\Psi_{w, \bar{\phi}}}^{y} D_{22} + B_{\Psi_{w, \bar{\phi}}}^{u}
	\end{bmatrix}, \nonumber \\
	\hat{C} = &
	\begin{bmatrix}
		D_{\Psi_{s,f}}^{y} C_1 & 0_{2d \times d} & 0_{2d \times d}\\
		 D_{\Psi_{w,f}}^{y} C_1 & C_{\Psi_{w,f}} & 0_{2d \times d}\\
		D_{\Psi_{s, \bar{\phi}}}^{y} C_2  & 0_{2d \times d} & 0_{2d \times d} \\
		D_{\Psi_{w, \bar{\phi}}}^{y} C_2  & 0_{2d \times d} & C_{\Psi_{w, \bar{\phi}}}
	\end{bmatrix},~ \nonumber\\
	\hat{D} = &
	\begin{bmatrix}
		D_{\Psi_{s, f}}^{y} D_{11} + D_{\Psi_{s, f}}^{u} & D_{\Psi_{s, f}}^{y} D_{12}\\
		D_{\Psi_{w, f}}^{y} D_{11} + D_{\Psi_{w, f}}^{u} & D_{\Psi_{w, f}}^{y} D_{12}\\
		D_{\Psi_{s, \bar{\phi}}}^{y} D_{21} & D_{\Psi_{s, \bar{\phi}}}^{y} D_{22} + D_{\Psi_{s, \bar{\phi}}}^{u} \\
		D_{\Psi_{w, \bar{\phi}}}^{y} D_{21} & D_{\Psi_{w, \bar{\phi}}}^{y} D_{22} + D_{\Psi_{w, \bar{\phi}}}^{u}
	\end{bmatrix},
\end{align}}
$B_1$, $B_2$, $C_1$, $C_2$, $D_{11}, \ldots, D_{22}$ are partitions of $B$, $C$, $D$ in \eqref{eq:system matrices discrete-time}, respectively;
The other matrices are given in \eqref{eq:weighted-off-by-one iqc state space}, \eqref{eq:sector iqc state space} with $K = L_f- \mu_f$, or $K = L_{\bar{\phi}} - \mu_{\bar{\phi}}$, depending on subscripts $f$, $\bar{\phi}$, respectively.
\end{theorem}
\hiro{The system matrices result from the standard factorization and state-space realization for the IQC condition \eqref{eq:IQC theorem condition 3} (\cite{scherer2000linear,seiler2015stability}).}
The proof is similar to \cite[Theorem 4]{lessard2016analysis} and is omitted here.
\revise{The convergence rate $\rho$ in \eqref{eq:LMI IQC discrete-time} needs to be treated as a constant such that 
\icrr{\eqref{eq:LMI IQC discrete-time}} is an LMI.  Then, a bisection search on $\rho$  can be carried out to obtain the optimal convergence rate.}

\subsection{Stepsize selection}\label{subsection Stepsize selection}
It is well-known that the optimal fixed stepsize for the GD method $x_{k+1} = x_{k} - \eta \nabla f (x_{k})$ is $\eta = \frac{2}{L_f + \mu_f}$, rendering the \li{tight} \ml{upper bound \iclr{for the} convergence rate} $\rho = \frac{\kappa_f - 1}{\kappa_f + 1}$, where $\kappa_f = L_f/\mu_f$ \revise{(\cite{nesterov2003introductory, ruszczynski2011nonlinear,lessard2016analysis})}.
\revise{However, the choice of constant stepsizes and \icrr{the} corresponding convergence rates for \icrr{the} MD method have not yet been \icrr{addressed} in the literature.}
\revise{\icrr{In our numerical studies in section \ref{Numerical Examples} we will use the optimal stepsize when the objective functions $f$, and the DGF $\phi$ are quadratic, since in this case the optimal stepsize and the corresponding tight upper bound on the rates can be analytically derived.
These are stated in the following proposition.}}
\revise{
\begin{proposition}\label{prop: MD quadratic convergence rate}
Let $f(x) = \frac{1}{2} x^T F x + p^T x + r$, and $\phi (x) = \frac{1}{2} x^T \Phi x$, where $p, r \in \mathbb{R}^{d}$ are constant vectors, $F$ and $\Phi$ are any positive definite matrices that satisfy $\mu_f I_d \leq F \leq L_f I_d$, and $\mu_{\phi} I_d \leq \Phi \leq  L_{\phi} I_d$, respectively.  Then, the smallest upper bound on the convergence rate for the MD method \eqref{eq:MD algorithm compact} is given by $\rho = \frac{\kappa_f \kappa_{\bar{\phi}} - 1}{ \kappa_f \kappa_{\bar{\phi}} + 1}$, where $\kappa_f=L_f/\mu_f$, $\kappa_{\bar{\phi}}=L_{\bar{\phi}}/\mu_{\bar{\phi}}$ with $\mu_{\bar{\phi}} = \left(L_{\phi} \right)^{-1}$, and $L_{\bar{\phi}} =\left( \mu_{\phi}\right)^{-1} $ are the condition numbers of $f$, $\bar{\phi}$, respectively, and this convergence rate is achieved with the stepsize $\eta = \frac{2}{L_f L_{\bar{\phi}} + \mu_f \mu_{\bar{\phi}}}$.
\end{proposition}
\begin{pf}
See Appendix~\ref{proof of prop: MD quadratic convergence rate}.
\end{pf}
}

\il{We will} \icrr{illustrate} 
via the numerical example in Section~\ref{Example: discrete-time} that \icrr{for} such a choice of $\eta$ \revise{
the associated tight convergence rate
bound $\rho = \frac{ \kappa_f \kappa_{\bar{\phi}} - 1}{\kappa_f \kappa_{\bar{\phi}} + 1}$
for quadratic functions \icrr{coincides with that obtained from the results in this paper} 
for any nonlinear functions $f \in S(\mu_f, L_f)$, and $\phi \in S(\mu_\phi, L_\phi)$. \icrr{\icrr{Also when different stepsizes were tried in our numerical computation of the rate bounds, we found that the smallest rate bound was obtained for this choice of stepsize.}}}

\icrr{A comparison} of the convergence rate bounds obtained from our results and those reported in the literature is also presented \icrr{in Section~\ref{Example: discrete-time}.}} 
\todoino{Is there a particular significance to the result described above? Is $\kappa$ the condition number of the composite nonlinearity - at least in certain cases? I think you mention in the examples section that this is a tight bound for quadratic functions?}

\section{\revise{Projected} MD and convergence rate}\label{Projected MD and convergence rate}
\lmm{The study on the \revise{projected} MD algorithm is more challenging due to the existence of \ic{a projection in the dynamics}. We \ic{will show in this section that the IQC framework developed in the previous section allows to derive bounds for the convergence rate in this case as well.}
}
\subsection{\revise{Projected} MD algorithm in the form of Lur'e systems}\label{sec:ProjMD}
The discrete-time MD algorithm \icl{with a convex} constraint set $\mathcal{X}$ is given by
\begin{align}\label{eq:MD algorithm set constraint compact}
x_{k+1} = \underset{x \in \mathcal{X}}{\text{argmin}} \left\{ \nabla f(x_k)^T x + \frac{1}{\eta} D_{\phi} (x, x_{k}) \right\}.
\end{align}

As commented \icl{in} \cite{beck2003mirror}, 
\icrr{one of the motivations} for the \icl{use of the Bregman divergence and the design of the distance function $\phi$} is to reflect the geometry of the given constraints set $\mathcal{X}$, \icrr{so
that the constraints can be satisfied without incorporating projections.
This is achieved when the function $\phi$ is 
chosen such that the following property holds}
\icrr{
\begin{equation}\label{eq:smooth assumption}
\begin{split}
	&\| \nabla \phi (x_k) \| \rightarrow +\infty \text{ as } k \rightarrow \infty, 
\\& \ \  \ \ \text{for all sequences } \  \{x_k \in \textrm{int}\mathcal{X}: \ k=1,2, \hdots\} \\ 
&\qquad \qquad \qquad \ \  \text{ where } x_k \rightarrow x \in \textrm{bd}\mathcal{X} \text{ as } k \rightarrow \infty.
\end{split}
\end{equation}
}
where $\textrm{bd}\mathcal{X}$ denotes the boundary of $\mathcal{X}$ \mmli{and $\textrm{int}\mathcal{X}$ represents the interior of $\mathcal{X}$}.
However, when \icrr{we use a function $\phi$ that does not satisfy}
assumption \eqref{eq:smooth assumption}, 
a \revise{projected} variant should be considered, \icrr{as in} 
\eqref{eq:MD algorithm set constraint compact} (\cite{beck2003mirror}).

\lmm{While \eqref{eq:MD algorithm compact} finds the optimal point by setting the gradient to zero, the optimization of \eqref{eq:MD algorithm set constraint compact} requires satisfying a differential inclusion to achieve optimality (\cite[Theorem 3.34]{ruszczynski2011nonlinear}), i.e., }
\begin{align}\label{eq:optimal solution of constrained MD}
	\revise{- \nabla \phi (x_{k+1}) + \nabla \phi (x_{k}) - \eta \nabla f(x_{k}) \in  N_{\mathcal{X}} (x_{k+1})}
\end{align}
where \lmm{$N_{\mathcal{X}} (x_{k+1})$ represents the normal cone of the convex set $\mathcal{X}$ at $x_{k+1}$.}
\lmm{From the strong convexity of $\phi$, we have that $x_{k+1}$ is the unique point satisfying the above inclusion \cite[Proposition 3.2]{beck2003mirror}.}
Then equivalently,
\begin{align}
	\nabla \phi (x_{k+1}) = \nabla \phi (x_{k}) - \eta \nabla f(x_{k}) - \nu T (x_{k+1})
\end{align}
where $\nu > 0$ is a positive parameter and \revise{$T (x)$ is an element in $N_{\mathcal{X}} (x)$ such that the above equation is satisfied}.
Again, letting $z_{k} = \nabla \phi (x_{k})$, we have that
\begin{align}\label{eq:constrained MD discrete-time composition}
	z_{k+1} = z_{k} - \eta \nabla f( \nabla \bar{\phi} (z_{k})) - \nu T \left( \nabla \bar{\phi} (z_{k+1}) \right)
\end{align}
where $\bar{\phi}$ is defined exactly the same as \eqref{eq:convex conjugate}, \ic{i.e.\hu{,} without} the set constraint.

\mm{Notice that there are two terms \ic{involving a} composition of nonlinear operators on the right hand side of \eqref{eq:constrained MD discrete-time composition}, which means that two transformations need to be carried out 
\ic{analogous to those in} the previous section. }
\lmm{Additionally, the presence of $z_{k+1}$ on the right-hand side necessitates an extra state to accommodate this algebraic loop.}
\ic{In order to incorporate the points above,} the \revise{projected} MD algorithm \eqref{eq:constrained MD discrete-time composition} is rewritten \ic{as} 
\begin{align}
	z_{k+1} = A z_{k} + B u_k, \quad y_k = C z_k + D u_k
\end{align}
with
\begin{align}\label{eq:system matrices discrete-time constrained MD}
& \begin{bmatrix}
	\begin{array}{c|c}
      A & B\\
      \hline
      C & D
    \end{array}
\end{bmatrix} \nonumber \\
=&
\begin{bmatrix}
	\begin{array}{c|cccc}
		(1-\eta \mu_f \mu_{\bar{\phi}})I_d &
		-\eta I_d & -\eta \mu_f I_d & - I_d & 0_d\\
		\hline
		\mu_{\bar{\phi}} I_d & 0_d & I_d & 0_d & 0_d\\
		I_d & 0_d & 0_d & 0_d & 0_d\\
	\mu_{\bar{\phi}} (1-\eta \mu_f \mu_{\bar{\phi}}) I_d & - \eta \mu_{\bar{\phi}} I_d & - \eta \mu_{\bar{\phi}} \mu_{f} I_d & -\mu_{\bar{\phi}} I_d & I_d \\
	(1-\eta \mu_f \mu_{\bar{\phi}})I_d & -\eta I_d & - \eta \mu_f I_d & - I_d & 0_d
	\end{array}
\end{bmatrix}
\end{align}
and system input
\begin{align}\label{eq:input discrete-time constrained MD}
	u_{k} = \begin{bmatrix}
u_k^{(1)} \\ u_k^{(2)} \\u_{k}^{(3)} \\ u_{k}^{(4)}
\end{bmatrix}
= &
\begin{bmatrix}
\nabla f (y_k^{(1)}) - \mu_{f} y_k^{(1)} \\ \nabla \bar{\phi} (y_k^{(2)}) - \mu_{\bar{\phi}} y_k^{(2)} \\ \nu T (y_k^{(3)}) \\
\nabla \bar{\phi} (y_k^{(4)}) - \mu_{\bar{\phi}} y_k^{(4)}
\end{bmatrix}.
\end{align}

Denote $z^\textup{\nkl{opt}}$ as the optimal \il{value of $z$ at} steady state, with corresponding \il{equilibrium values} $y^\textup{\nkl{opt}}$, and $u^\textup{\nkl{opt}}$. It is worth noting that \lmm{in this reformulation}, $y_{k}^{(1)} = x_{k}$, $y_{k}^{(2)} = z_{k}$, $y_{k}^{(3)} = z_{k+1}$ and $y_{k}^{(4)} = x_{k+1}$. Then, we have $y^{(2), \textup{opt}} = y^{(4), \textup{opt}}$.

Define $\tilde{u}_{k} = u_k - u^\textup{\nkl{opt}}$ and $\tilde{y}_{k} = y_k - y^{\textup{opt}}$.
Obviously, the input/output pair $(\tilde{u}_k$, $\tilde{y}_{k})$
of the nonlinearity satisfies the Sector IQC and weighted off-by-one IQC with suitable parameters. \hiro{In particular, \ir{the normal cone of a convex set $\mathcal{X}$ at a point $x\in\mathcal{X}$, is equal to 
the subdifferential of the indicator function on $\mathcal{X}$, which is 
convex} \cite[p.~215]{rockafellar1997convex}. Thus, we have
\begin{align*}
	\left( T (x) - T(y) \right)^T (x - y) \geq 0,
 \ir{\forall x\in\mathcal{X}, \forall y\in\mathcal{X}},
\end{align*}}
\ir{which is a property that also follows from the definition of the normal cone.}
Then, $T(x)$ satisfies \ir{for all trajectories in $\mathcal{X}$}
the Sector and weighted off-by-one IQCs with $m = 0$ and $L = \infty$ (\cite{lessard2016analysis}).
We include the state-space realization for $\Psi_{s, T}$, $\Psi_{w, T}$ as follows,
\begin{align*}
\Psi_{s,T}:
	 & \begin{bmatrix}
		\begin{array}{c|cc}
		A_{\Psi_{s, T}} & B_{\Psi_{s, T}}^{y} & B_{\Psi_{s, T}}^{u}\\
		\hline
		C_{\Psi_{s, T}} & D_{\Psi_{s, T}}^{y} & D_{\Psi_{s, T}}^{u}
		\end{array}
	\end{bmatrix}
	\ic{:=} \begin{bmatrix}
		\begin{array}{c|cc}
		0_d & 0_d & 0_d\\
		\hline
		0_d & 0_d & I_d\\
		0_d & I_d & 0_d
		\end{array}
	\end{bmatrix},\\
\Psi_{w, T}:
	 & \begin{bmatrix}
		\begin{array}{c|cc}
		A_{\Psi_{w, T}} & B_{\Psi_{w, T}}^{y} & B_{\Psi_{w, T}}^{u}\\
		\hline
		C_{\Psi_{w, T}} & D_{\Psi_{w, T}}^{y} & D_{\Psi_{w, T}}^{u}
		\end{array}
	\end{bmatrix}
	\ic{:=} \begin{bmatrix}
		\begin{array}{c|cc}
		0_d & -I_d & 0_d\\
		\hline
		0_d & 0_d & I_d\\
		\bar{\rho}^2 I_d & I_d & 0_d
		\end{array}
	\end{bmatrix}.
\end{align*}
The above weighted off-by-one IQC in the frequency domain is given by \cite{boczar2017exponential}.

\subsection{Repeated nonlinearities}
We observe that there are repeated nonlinearities present in equation \eqref{eq:input discrete-time constrained MD}.
To obtain less conservative results, these repeated nonlinearities must be properly addressed.
The largest class of bounded \lmm{convolution} multipliers to deal with repeated nonlinearities
have been proposed in the literature (\cite{kulkarni2002all,d2001new}). \li{In this paper, we opt for a special class of IQCs described below, which are sufficient for providing feasible solutions with a good convergence rate bound.}
\todoino{the previous sentence does not read very well}

Notably, $y_k^{(2)} - y_k^{(4)}$ and $u_k^{(2)} - u_k^{(4)}$, \lmm{as the input and output of a nonlinearity}, also satisfy the Sector IQC in Lemma~\ref{lem:co-coercivity sector IQC}. It is evident that they also satisfy the weighted-off-by-one IQC in Lemma~\ref{lem:weighted-off-by-one IQC}.
As a result, we obtain two additional IQCs necessary for describing the repeated nonlinearities in \eqref{eq:input discrete-time constrained MD}.

\li{It is worth mentioning that it may be possible to select additional classes of IQCs from \cite{kulkarni2002all,d2001new} in order to achieve less conservative rates. However, exploring these possibilities falls beyond the scope of this paper.}


\subsection{Convergence rate analysis}
\revise{
Building upon previous discussions,  for each nonlinearity in \eqref{eq:input discrete-time constrained MD}, we incorporate one Sector IQC and one weighted off-by-one IQC, yielding $ 2 \times 4 = 8$ IQCs.
Additionally, for repeated nonlinearities, we also use both the Sector and weighted off-by-one IQCs, amounting to $2 \times 1 = 2$ IQCs.
Altogether, this results in a conic combination of $10$ IQCs to characterize the nonlinearity in  \eqref{eq:input discrete-time constrained MD}. A theorem stemming from this discussion, analogous to Theorem~\ref{thm:convergence rate discrete-time}, can be obtained with its proof omitted for brevity.}
\begin{theorem}\label{thm:constrained MD convergence rate}
The mirror descent algorithm converges linearly with rate \hiro{$\bar{\rho} \leq \rho \leq 1$} if there exists a matrix \hiro{$P = P^T > 0$}, $\alpha = \textup{diag} \{\alpha_1,\ldots,\alpha_{10}\} \geq 0$ such that the following inequality holds
	\begin{align}\label{eq:LMI IQC discrete-time constrained MD}
	\left[
	\begin{smallmatrix}
		\hat{A}^T P \hat{A} - \rho^2 P & \hat{A}^T P \hat{B}\\
		* & \hat{B}^T P \hat{B}
	\end{smallmatrix} \right]
	+
	\left[\begin{smallmatrix}
		\hat{C} & \hat{D}
	\end{smallmatrix} \right]^T
	\left( \alpha \otimes M \right)
	\left[
	\begin{smallmatrix}
		\hat{C} & \hat{D}
	\end{smallmatrix} \right]
	\leq 0,
	\end{align}
where $\hat{A}$, $\hat{B}$, $\hat{C}$, $\hat{D}$ are given in Appendix~\ref{Matrices for LMI}.
\end{theorem}
\hiro{Similar to Theorem~\ref{thm:convergence rate discrete-time}, the system matrices are obtained through the standard factorization and state-space realization for the IQC condition \eqref{eq:IQC theorem condition 3}.}
\lmm{To the best of our knowledge, this theorem
provides \ic{the first systematic method for computing 
convergence rate bounds for the} \revise{projected} MD method \eqref{eq:MD algorithm set constraint compact} 
\ic{with} a constant stepsize. These bounds are further illustrated in the next section via numerical examples.}

\section{Numerical Examples}\label{Numerical Examples}
In this section, we present \revise{four} numerical examples to illustrate the IQC analysis for the continuous-time, discrete-time and \revise{projected} MD algorithm, respectively.
\subsection{Continuous-time MD method}\label{Example: Continuous-time}
\todoing{\st{what do you mean by feasibility ranges?}}
\ml{We investigate and compare the \il{feasibility}} of the IQC condition \eqref{eq:IQC theorem condition 3} \il{when} using merely the Sector IQC defined by \eqref{eq:IQC sector bounded}, and when using the conic combination of the Sector and Popov IQCs defined by \eqref{eq:IQC sector bounded + Popov}.
\ml{The frequency-\icl{domain} condition \eqref{eq:IQC theorem condition 3} under \eqref{eq:IQC sector bounded} can be easily transformed into \li{an LMI}
via the Kalman–Yakubovich–Popov (KYP) lemma. 
\icc{Condition} \eqref{eq:IQC theorem condition 3} under \eqref{eq:IQC sector bounded + Popov} is satisfied if and only if \eqref{eq:LMI exponential rate continuous-time} in Theorem~\ref{thm:convergence rate continuous-time} is feasible for some $\rho > 0$.}
Let $\eta = 1$, $\mu_f = 1$ and $\mu_{\bar{\phi}} = 1$, and $L_f = L_{\bar{\phi}}$.
\il{The feasibility of the IQCs (for some $\rho>0$) with varying composite condition number $\kappa = \frac{L_f}{\mu_{f} } \cdot \frac{L_{\bar{\phi}}}{\mu_{\bar{\phi}}}$ is} shown in Fig.~\ref{fig:continuous-time feasibility}.
\icrr{It should be noted that the largest convergence rate bound is $\rho = \eta \mu_{f} \mu_{\bar{\phi}}$ (Theorem \ref{thm:convergence rate continuous-time}) and is obtained with a conic combination of \mmli{Sector and Popov IQCs}.}}
Note that the MD algorithm should converge for any $L_f > \mu_f$ and $L_{\bar{\phi}} > \mu_{\bar{\phi}}$. However, we can observe that the Sector IQC defined by \ml{\eqref{eq:IQC sector bounded}} fails to certify the convergence of the MD method for $\kappa \geq 34$.
On the other hand, using the conic combination of the Sector and Popov IQCs defined by \ml{\eqref{eq:IQC sector bounded + Popov}} suffices to certify its convergence for arbitrary $\kappa$.

\mm{Since the Popov IQC corresponds to the Bregman divergence function in time-domain analysis,
this example thus \lmm{highlights} the importance of adopting Bregman-like Lyapunov functions in the convergence analysis.}

\begin{figure}[htbp]
	\centering
		\includegraphics[width = 1\linewidth]{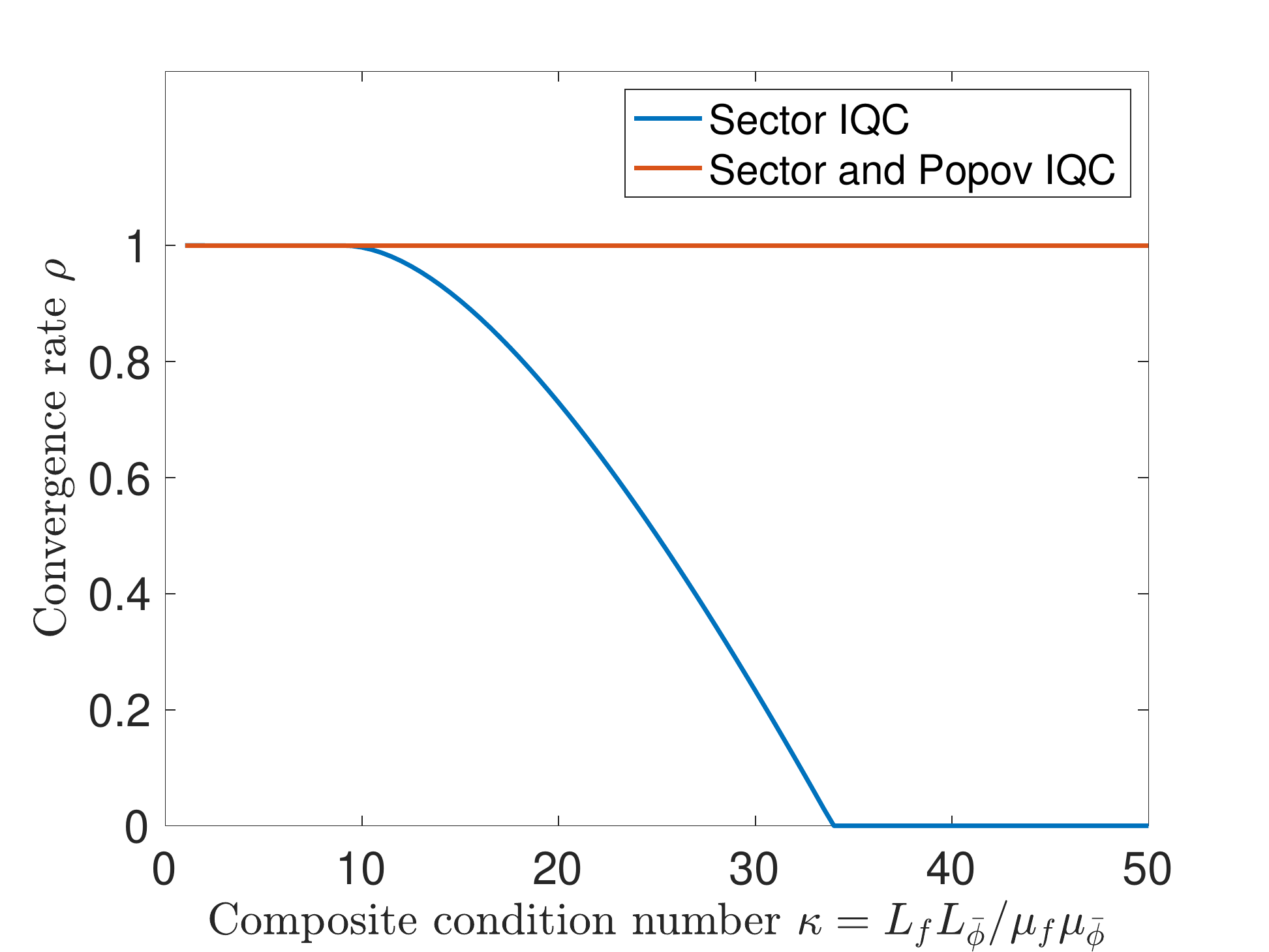}
	\caption{\hiro{Convergence rate for the continuous-time MD method as a function of
the composite condition number $\kappa =  \frac{L_f}{\mu_{f} } \cdot \frac{L_{\bar{\phi}}}{\mu_{\bar{\phi}}}$, derived from the Sector IQC \eqref{eq:IQC sector bounded} and Sector and Popov IQC \eqref{eq:IQC sector bounded + Popov}, respectively}.
\irr{The MD method parameters have been set as  $\eta =\mu_f=\mu_{\bar{\phi}} = 1$.}
\hiro{\icrr{Note that the largest convergence rate bound is $\rho = \eta \mu_f \mu_{\bar{\phi}} = 1$, as established in Theorem~\ref{thm:convergence rate continuous-time}, and we let $\rho = 0$ when the associated LMI is infeasible.}}
}
	\label{fig:continuous-time feasibility}
\end{figure}

\subsection{Discrete-time MD method}\label{Example: discrete-time}
Next, we present the convergence rate for the discrete-time MD method.
Let $\mu_f = 1$, $\mu_{\bar{\phi}} = 1$, and $L_f = L_{\bar{\phi}}$. Let the stepsize be $\eta = \frac{2}{L_f L_{\bar{\phi}} + \mu_f \mu_{\bar{\phi}}}$ as Section~\ref{subsection Stepsize selection} suggested. We compare the smallest convergence rate obtained from \eqref{eq:LMI IQC discrete-time} in Theorem~\ref{thm:convergence rate discrete-time} with that obtained from the semidefinite programs (SDPs) in \cite{sun2022centralized}, where the stepsize and convergence rate are both decision variables. The SDPs in \cite{sun2022centralized} are derived from the Lyapunov function $V (z_k) = \rho^{-k} D_{\bar{\phi}}(z_{k}, z^\textup{\nkl{opt}})$, \todoinc{is $\rho$ in the definition of $V$ the convergence rate - if not use a different symbol?\\ Yes it is the convergence rate.} which \il{is the} Bregman divergence function when $\rho = 1$.
The relation between the composite condition number $\kappa = \frac{L_f}{\mu_{f} } \cdot \frac{L_{\bar{\phi}}}{\mu_{\bar{\phi}}}$ and the convergence rate $\rho$ is shown in Fig.~\ref{fig:discrete-time convergence rate}.
We can observe that using the IQC analysis provides a tighter bound \il{for the} convergence rate.
\icc{We remark that the convergence rate bound obtained in this example using the IQC analysis coincides with the curve $\rho =\frac{ \kappa - 1}{\kappa  + 1 }$. This illustrates that this bound is also tight since \revise{Proposition~\ref{prop: MD quadratic convergence rate} shows that} this expression is the smallest upper bound on the convergence rate}
for all quadratic functions $f \in S(\mu_f, L_f)$ and $\bar{\phi} \in S(\mu_{\bar{\phi}}, L_{\bar{\phi}})$}.
\begin{figure}[htbp]
	\centering
	\includegraphics[width = 1\linewidth]{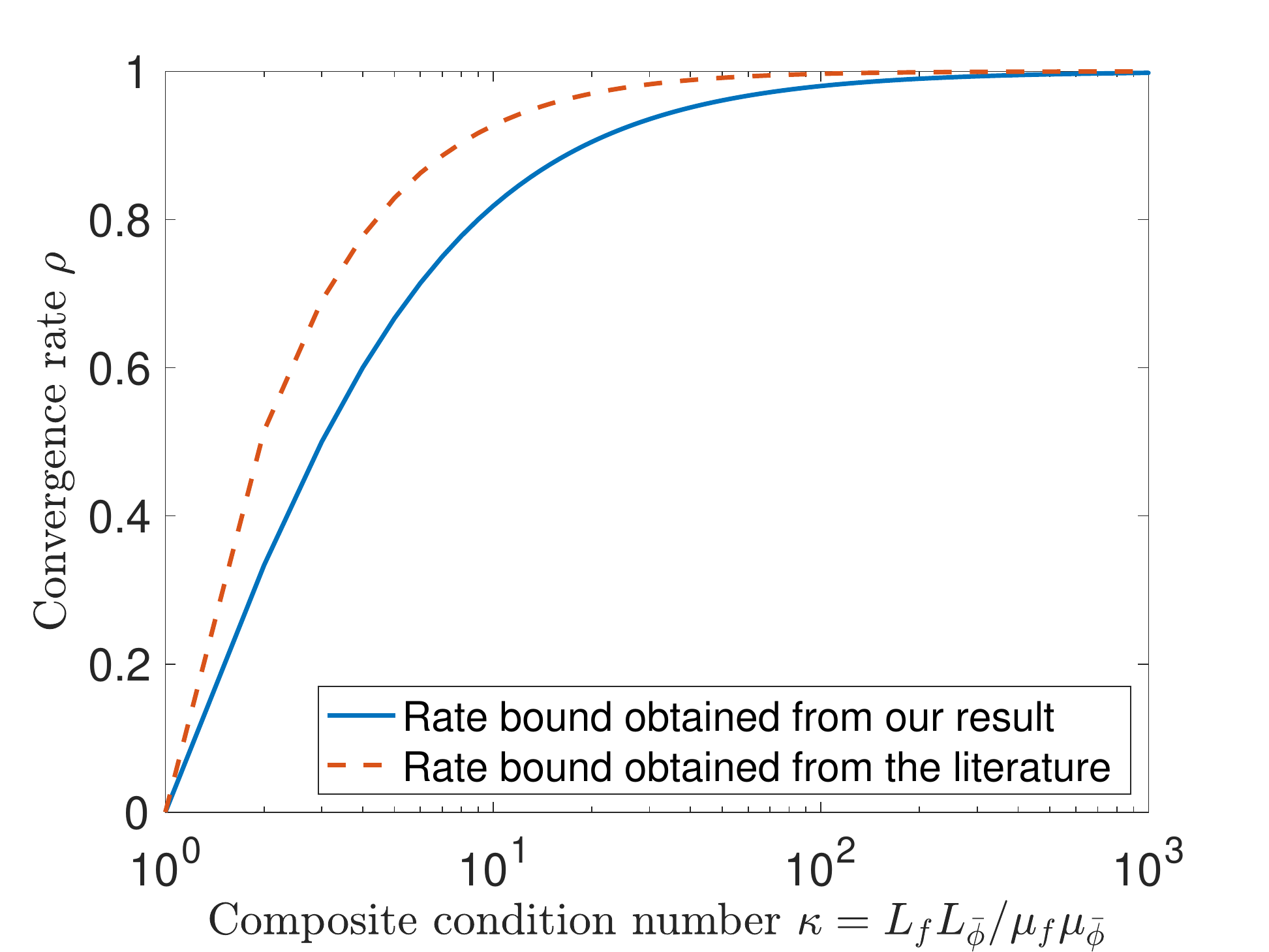}
	\caption{Convergence rate obtained from \eqref{eq:LMI IQC discrete-time} in Theorem~\ref{thm:convergence rate discrete-time} and from the SDPs in \cite{sun2022centralized}.  The convergence rate bound obtained from our result coincides with the curve $\rho =\frac{ \kappa - 1}{\kappa  + 1 }$.}
	\label{fig:discrete-time convergence rate}
\end{figure}
\revise{
\subsection{Convergence rate with quadratic functions}}
\mm{When comparing the convergence rate bounds between the GD and MD algorithms, we have $\frac{\kappa_f - 1 }{ \kappa_f + 1} \leq \frac{\kappa -1 }{ \kappa + 1}$. This is due to $\kappa = \kappa_f \kappa_{\bar{\phi}} \geq \kappa_f$. Hence, the distance generating function $\phi$ used in the MD method does not necessarily result in a faster convergence rate than the GD method under constant stepsizes.}
\mml{
The potential for improvement in the convergence rate of the MD algorithm is dependent on the specific structure of the distance generating function $\phi$, which must be carefully tailored to the function $f$.

\revise{
To illustrate this, we present a simple example that compares the convergence rate of the GD and MD methods when applied to a quadratic convex objective function. Specifically, we consider the following optimization problem
\begin{align*}
\underset{x \in \mathbb{R}^2}{\min} \left\{ f(x) := \frac{1}{2} x^T F x + p^T x
\right\}
\end{align*}
where
\begin{align*}
F =
\begin{bmatrix}
100 & -1\\ -1 & 1
\end{bmatrix}, \quad
p =
\begin{bmatrix}
1 \\ 10
\end{bmatrix}.
\end{align*}
The optimal solution and value can be  easily obtained,
\begin{align*}
x^{\textup{opt}} = - F^{-1} p =
\begin{bmatrix}
-0.1111 \\ -10.1111
\end{bmatrix},
 \quad f^{\textup{opt}} = -50.6111.
\end{align*}
It is clear that $f \in S(\mu_f, L_f)$ with $\mu_f = 0.9899$ and $L_f = 100.0101$.
Next, we select DGF $\phi (x) = \frac{1}{2} x^T \Phi x$, where
$
\Phi =
\begin{bmatrix}
10 & 1\\
1 & 1
\end{bmatrix}
$. As a result, its convex conjugate is $\bar{\phi} (z)= \frac{1}{2}  z^T \bar{\Phi} z$, and $\bar{\phi} (z) \in S(\mu_{\bar{\phi}}, L_{\bar{\phi}})$, where $\mu_{\bar{\phi}} = 0.09891$ and $L_{\bar{\phi}} =  1.1233$.
The convergence rate bound of the GD algorithm is $\rho_g = \frac{L_f/\mu_f - 1}{L_f/\mu_f + 1} = 0.9804$, while the convergence rate bound of the MD algorithm given by Proposition~\ref{prop: MD quadratic convergence rate} (which holds for arbitrary quadratic functions in the specified classes) is
$$\rho = \frac{L_f L_{\bar{\phi}}/ \mu_f \mu_{\bar{\phi}} - 1}{L_f L_{\bar{\phi}}/ \mu_f \mu_{\bar{\phi}} + 1} = 0.9982 > \rho_g.$$
Interestingly, due to the specific structures of $F$ and $\bar{\Phi}$, the spectrum for $F\bar{\Phi}$ in \eqref{eq: spectrum bound} becomes $0.9576 \leq \lambda \leq 11.4868$. As such, a more precise convergence rate bound is $$\rho_{m} = \frac{11.4868/0.9576 - 1}{11.4868/0.9576 + 1} = 0.8461 < \rho_g,$$ indicating an accelerated convergence compared to the GD method.
This observation is supported by the optimization errors $|f(x_k) - f^{\textup{opt}}|$ for both the GD and MD methods, depicted in Fig.~\ref{fig:QP rate comparison}.
Applying the MD method here can be viewed as a form of the preconditioning gradient descent (\cite{maddison2021dual}). For more general objective functions, investigating the potential for accelerated convergence rates using the MD method will be an interesting extension for future research.
\begin{figure}[htbp]
\centering
\includegraphics[width = 1\linewidth]{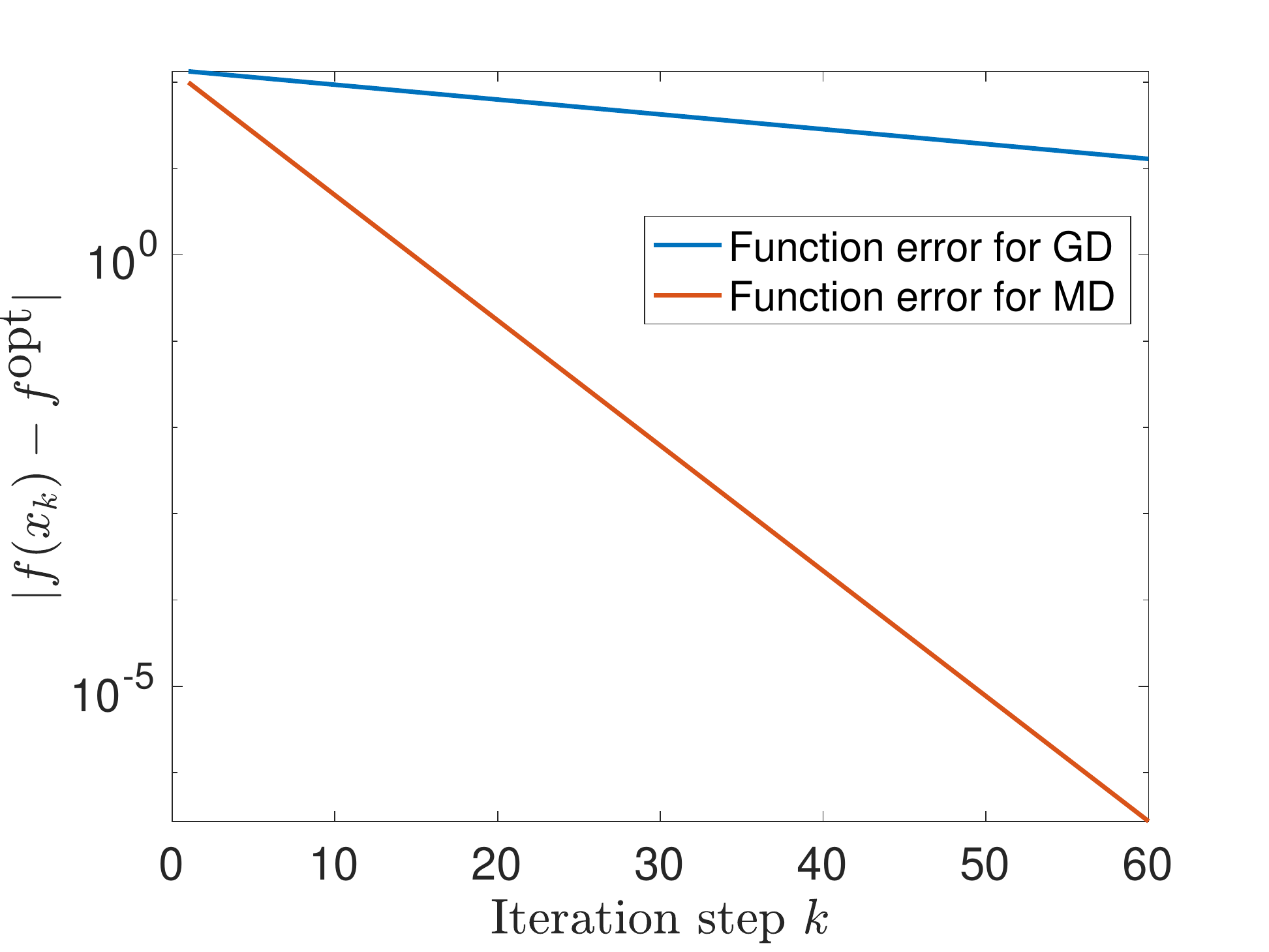}
\caption{Optimization errors for GD and MD methods, respectively.}
\label{fig:QP rate comparison}
\end{figure}
}

\subsection{\revise{Projected} MD method}
Lastly, we present the convergence rate of the \revise{projected} MD algorithm.
Let $L_{f} = L_{\bar{\phi}} = 1$ and $\mu_{f} = \mu_{\bar{\phi}}$.
We demonstrate the relation between the composite condition number $\kappa = \frac{L_f}{\mu_{f} } \cdot \frac{L_{\bar{\phi}}}{\mu_{\bar{\phi}}}$ and convergence rate $\rho$ using Theorem~\ref{thm:constrained MD convergence rate}.
\revise{For the stepsize selection, we start by looking
for the
\icrr{closest 
stepsize to $\eta = \frac{2}{L_f L_{\bar{\phi}} + \mu_f \mu_{\bar{\phi}}}$ for which the LMIs are feasible,}
 through a bisection search with fixed $\rho = 1$. After finding a suitable stepsize, we then use the bisection method again to search for
\icrr{smallest rate bound.}}
The result is depicted in Fig.~\ref{fig:constrained MD}, \hiro{showing the smallest rate bounds we can obtain from the proposed IQCs, \irr{after a search over different stepsizes was also carried out}.}
\begin{figure}[htbp]
	\centering
	\includegraphics[width = 1\linewidth]{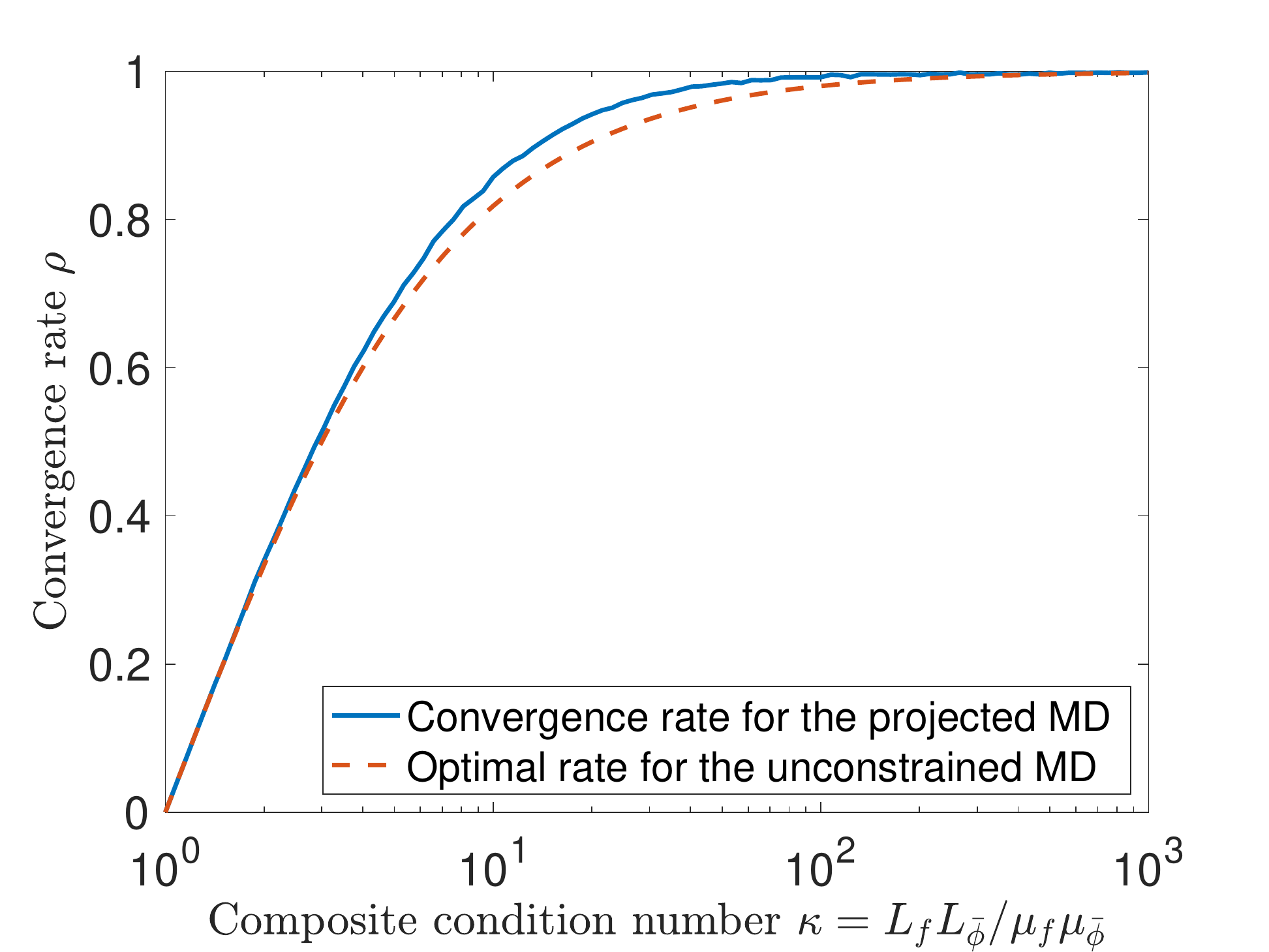}
	\caption{Convergence rate for the \revise{projected} MD algorithm \eqref{eq:MD algorithm set constraint compact} obtained from Theorem~\ref{thm:constrained MD convergence rate}.}
	\label{fig:constrained MD}
\end{figure}
\lmm{It is important to note that the convergence rate we obtain in this case is \li{slightly larger than the one in the} unconstrained scenario.
This is due to the presence of additional nonlinear term $\nu T (x_{k+1})$ in \eqref{eq:constrained MD discrete-time composition} and the repeated nonlinearities in \eqref{eq:input discrete-time constrained MD}. As a direction for future work, it would be intriguing to investigate the possibility of obtaining
\ic{a tighter}
convergence rate bound by employing more advanced IQCs.}

\section{Conclusion}\label{Conclusion}
An integral quadratic constraint analysis framework has been established for the mirror descent method in both continuous-time and discrete-time settings.
 \lmm{In the continuous-time setting,} we demonstrated that the Bregman divergence function is a special case of the \icl{Lyapunov functions associated with the Popov criterion \ill{when these are applied to an appropriate reformulation of the \revise{MD dynamics}}.}
\lmm{In the discrete-time setting,} \ml{we established \ill{upper bounds} for
\ic{the 
convergence}
\ill{rate via appropriate IQCs applied to the transformed system. It has also been illustrated via numerical examples that the convergence rate bound for the unconstrained case is tight.}}
Future endeavors \il{include extending the framework developed to other related algorithms such as accelerated MD methods.}
\todoing{\st{not sure what you mean by the previous phrase}\\
Answer: The optimal convergence rate is only provided numerically. I would like to show that it is such analytically.}




\begin{ack}                               
This work was partly supported by ERC starting grant 679774, and \ic{the} Japanese Ministry of the Environment.
M.~Li is also supported by Japan Society for the Promotion of Science (JSPS) KAKENHI under Grant 24K23864.
For the purpose of open access, the authors have applied a Creative Commons Attribution (CC BY) licence  to any Author Accepted Manuscript version arising.
\end{ack}

\bibliography{References}           



\section*{Appendix}
\appendix
\section{Proof of Theorem~\ref{thm:convergence continuous-time IQC}}\label{Proof of thm:convergence continuous-time IQC}    
We prove this theorem by applying \irr{the extension of Theorem~\ref{thm:IQC theorem} in
\cite[Theorem 1]{jonsson1997stability} that allows the consideration of conic combinations of IQCs with essentially bounded multipliers $\Pi(j\omega)$, \ir{with} the Popov IQC.}
\ml{\icl{The} closed-loop system with $\tau \Delta$ is well-posed since the system we are investigating is the MD algorithm with $\tau$ interpreted as \li{a scaling of the} gradients.}
Then, the first condition of Theorem~\ref{thm:IQC theorem} is satisfied.
We use a conic combination the Sector and Popov IQCs \kl{defined by $\Pi$ and $\Pi_P(j\omega)$ given by~\eqref{eq:IQC sector bounded} and~\eqref{eq:Popov IQC for delta}, respectively}.
\icl{In particular,} $\Delta \in \text{IQC} \left( \Pi \right)$, where $\Pi ( j \omega )$ takes the form
\begin{align}\label{eq:IQC sector bounded + Popov}
	\Pi ( j \omega ) = &
	\left[
	\begin{smallmatrix}
	0_d & 0_d & \left( \alpha_1 ( L_f - \mu_f ) - \beta_1 j \omega \right) I_d & 0_d \\
	0_d & 0_d & 0_d & \left( \alpha_2 ( L_{\bar{\phi}} - \mu_{\bar{\phi}} ) - \beta_2 j \omega \right) I_d\\
	* & * & -2 \alpha_1 I_d & 0_d \\
	* & * & 0_d & -2 \alpha_2 I_d
\end{smallmatrix}
\right] \nonumber \\
:= & \begin{bmatrix} \Pi_{11} ( j \omega ) & \Pi_{12} ( j \omega ) \\ \Pi_{12}^* ( j \omega ) & \Pi_{22} ( j \omega ) \end{bmatrix}.
\end{align}
Since $\Pi_{11} ( j \omega ) \geq 0$ and $\Pi_{22} ( j \omega ) \leq 0$, the second condition of Theorem~\ref{thm:IQC theorem} is satisfied.
Condition \eqref{eq:IQC theorem condition 3} \irr{in Theorem~\ref{thm:IQC theorem}} is satisfied if
\begin{equation}\label{eq:iqc inequality}
\begin{aligned}
& -\beta_1^2 \omega^4 +
\left( - \eta^2 \beta_2^2 + 2\alpha_1 \eta( L_f + \mu_f) \beta_2 \right.\\
& \left. + 4 \alpha_1 \alpha_2 - \alpha_1^2 (L_f - \mu_f)^2 + 2  \alpha_2 \beta_1 \eta L_{\bar{\phi}}  + 2 \alpha_2 \beta_1 \eta \mu_{\bar{\phi}} \right) \omega^2\\
& + \eta^2 \left( 4 \alpha_1  \alpha_2 L_f \mu_f L_{\bar{\phi}} \mu_{\bar{\phi}} - \alpha_2^2 (L_{\bar{\phi}} - \mu_{\bar{\phi}} )^2 \right) > 0, ~ \forall \omega \in \mathbb{R}.
\end{aligned}
\end{equation}
Clearly, a necessary condition for the above inequality to hold is $\beta_1 = 0$, \revise{which also coincides with the multiplier condition in \cite{jonsson1997stability}}. To show that it is feasible, we let $\alpha_2 = 1$,  $\alpha_1$, $\beta_2$ take the values of $\alpha_1^* = \frac{(L_{\bar{\phi}} - \mu_{\bar{\phi}})^2}{4 L_f L_{\bar{\phi}} \mu_f \mu_{\bar{\phi}}}$, and $\beta_2^* = \frac{\alpha_1^*(L_f + \mu_f) + 2 \sqrt{ {\alpha_1^*}^2 L_f \mu_f + \alpha_1^*}}{\eta}$, respectively, and the left hand side of \eqref{eq:iqc inequality} becomes zero. Then, there exist $0 < \beta_2 < \beta_2^*$ and $\alpha_1 > \alpha_1^*$ such that \eqref{eq:iqc inequality} is satisfied for all $\omega \in \mathbb{R}$.
The stability implies $v \in \mathbf{L}_2[0, \infty)$ for $e, g \in \mathbf{L}_2 [0, \infty)$ in \eqref{eq:feedback interconnection model}.
\mm{For the error system \eqref{eq:mirror descent algorithm continuous-time}, the initial condition response of $z_0 - z^\textup{\nkl{opt}}$ can be represented by the external input $g = C e^{A t} \left( z_0 - z^\textup{\nkl{opt}}\right)$, which satisfies $g \in \mathbf{L}_2[0, \infty)$ and $g \to 0$ as $t \rightarrow \infty$.  This means that $v \rightarrow 0$ as $t \rightarrow \infty$ according to \cite[Proposition 1]{megretski1997system},  and the \lmm{system trajectories converge to zero.}
Thus, the trajectory of $x = \nabla \bar{\phi} (z)$ for \il{any input} $z_0$ converges to the optimal solution of problem \eqref{eq:optimization problem}.}
\todoing{\st{IQCs are about input/ouput stability. The fact that this implies convergence for any initial condition needs some further considerations, see e.g. Megretski-Rantzer paper.}}

\section{Proof of Theorem~\ref{thm:convergence rate continuous-time}}\label{Proof of thm:convergence rate continuous-time}
Consider \ir{in \eqref{eq:mirror descent algorithm continuous-time}, \eqref{eq:nonlinearity} the transformation $\tilde{z}_{\rho}(t):=e^{\rho t}\tilde z(t)$,  the modified transfer function} $G_{\rho} (s): = C (s I - A - \rho I)^{-1} B + D$ with $(A, B, C,D)$ given in \eqref{eq:transfer function continuous-time}, and $\Delta_{\rho}(t, v(t)) := e^{\rho t} \Delta (e^{-\rho t} v(t))$ \hirot{with $\Delta$ given in \eqref{eq:nonlinearity}}. \ir{Note that  
$G_{\rho} (s)$ is the transfer function of the system with input  $\tilde u_\rho(t):= e^{\rho t} \hiro{\tilde{u}}(t)$, output $\tilde y_\rho(t):= e^{\rho t} \hiro{\tilde{y}} (t)$, and state vector~$\tilde{z}_{\rho}(t)$.}

\ir{According to \cite{hu2016exponential}, the closed-loop system of $G(s)$ and $\Delta$ \ir{in \eqref{eq:mirror descent algorithm continuous-time}, \eqref{eq:nonlinearity}} is exponentially stable with rate $\rho$ \hiro{iff} the closed-loop system of $G_{\rho} (s)$ and $\Delta_{\rho}$
is \hirot{linearly} stable, \hirot{i.e., $\| \tilde{z}_{\rho} (t) \| \leq c \| \tilde{z}_{\rho} (0)\|$, $\forall t \geq 0$, \ir{for some $c > 0$.}}}

\ir{In the remainder of the proof we show via a Lyapunov function that this linear stability property holds when the LMI \eqref{eq:LMI exponential rate continuous-time} holds.}
\hirot{
Consider the \ir{following Lyapunov-like} function 
\ir{for the transformed system}
\begin{align}\label{eq:lyapunov function popov 0}
V_{\rho} = \frac{1}{2} \tilde{z}_{\rho}^T P \tilde{z}_{\rho} + \gamma \int_{0}^{t}  e^{\rho \tau } \Delta_2 ({ e^{-\rho \tau } \tilde{y}^{(2)}_{\rho} (\tau) } )^T \dot{\tilde{y}}^{(2)}_{\rho}(\tau)  d \tau
\end{align}
where $P = P^T > 0$ satisfying \eqref{eq:LMI exponential rate continuous-time}.\\
\hiro{The integral term can be rewritten as the line integral
$\gamma \int_{0}^{{\tilde{y}}^{(2)}_{\rho}(t)}  \Delta_{\rho}^{(2)} (\tau, \tilde{y}_{\rho} (\tau)) d {\tilde{y}}^{(2)}_{\rho}(\tau)$.
As $\Delta_{\rho}^{(2)}$ is slope-restricted with respect to ${\tilde{y}}_{\rho}^{(2)}(t)$ on sector $[0, L_{\bar{\phi}} - \mu_{\bar{\phi}}]$,  it \ir{follows
that}
\begin{align*}
0 \leq \gamma \int_{0}^{{\tilde{y}}^{(2)}_{\rho}(t)} \Delta^{(2)}_{\rho} (\tau, \tilde{y}_{\rho} (\tau) ) d {\tilde{y}}^{(2)}_{\rho}(\tau)\\
\leq \frac{ \gamma\left( L_{\bar{\phi}} - \mu_{\bar{\phi}}\right) }{2} \| {\tilde{y}}^{(2)}_{\rho}(t) \|^2.
\end{align*}
Moreover,
\lir{due to the structure of $D$ in \eqref{eq:system matrices continuous-time}, $(I - D \Delta_{\rho})$ is invertible with the inverse having a bounded gain, hence}
$(I - D \Delta_{\rho} ) \tilde{y}_{\rho} = C \tilde{z}_{\rho}$ implies $\| \tilde{y}_{\rho} \| \leq  \left\| (I - D \Delta_{\rho} )^{-1} C \right\|  \cdot \| \tilde{z}_{\rho} \|$.}
Thus, 
\ir{
we have the following linear bounds \lir{on $V_{\rho}$ } w.r.t. the Euclidean norm of the state vector $\tilde{z}_{\rho}$, }
 \hiro{
$
\alpha_{1} \| \tilde{z}_{\rho} \|^2 \leq V_{\rho} \leq \alpha_{2} \| \tilde{z}_{\rho} \|^2
$
for some $\alpha_1$, $\alpha_2 > 0$. \lir{We also have}}
\begin{align*}
	\dot{V}_{\rho} = & \tilde{z}_{\rho}^T P ( \tilde{A}_{\rho} \tilde{z}_{\rho} - \tilde{B} \tilde{u}_{\rho} ) + \left( \Delta_{\rho}(t,  \tilde{y}_{\rho}) \right)^T \Gamma C ( \tilde{A}_{\rho} \tilde{z}_{\rho} - \tilde{B} \tilde{u}_{\rho}  ) \\
	= & \frac{1}{2}
	\begin{bmatrix}
		\tilde{z}_{\rho} \\ \mm{-} \tilde{u}_{\rho}
	\end{bmatrix}^T
	\begin{bmatrix}
		P \tilde{A}_{\rho} + \tilde{A}_{\rho}^T P & P\tilde{B} - \tilde{C}_{\rho}^T \\ * & -(\tilde{D} + \tilde{D}^T)
	\end{bmatrix}
	\begin{bmatrix}
		\tilde{z}_{\rho} \\ - \tilde{u}_{\rho}
	\end{bmatrix}  \\
	& \mm{-} \tilde{z}_{\rho}^T \tilde{C}_{\rho}^T \tilde{u}_{\rho} + \tilde{u}_{\rho}^T \tilde{D} \tilde{u}_{\rho} + \tilde{u}_{\rho}^T \Gamma C ( \tilde{A}_{\rho} \tilde{z}_{\rho} \mm{-} \tilde{B} \tilde{u}_{\rho} )\\
	\leq & -  \tilde{u}_{\rho}^T \left( \mm{( Q C + \rho \Gamma C )} \tilde{z}_{\rho}  -  Q D  \tilde{u}_{\rho}+ Q  K^{-1}  \tilde{u}_{\rho} \right) \\
	= & -  \tilde{u}_{\rho}^T Q \left( \tilde{y}_{\rho} - K^{-1}  \tilde{u}_{\rho} \right) \mm{- \rho  \tilde{u}_{\rho}^T \Gamma C \tilde{z}_{\rho} } \\
	= & -  \tilde{u}_{\rho}^T Q \left( \tilde{y}_{\rho} - K^{-1}  \tilde{u}_{\rho} \right) - \rho  \tilde{u}_{\rho}^T \Gamma \left ( \tilde{y}_{\rho} - D \tilde{u}_{\rho} \right) \\
	\leq & 0
\end{align*}
where $\tilde{u}_{\rho} = \Delta_{\rho} (t, \tilde{y}_{\rho}) = e^{\rho t} \Delta(\tilde{y}_{\rho} e^{-\rho t})$ is the input of $G_{\rho} (s)$,
the first equality follows from $\Gamma C = \begin{bmatrix}
	0_d \\ \gamma I_d
\end{bmatrix}$, the first inequality follows from \eqref{eq:LMI exponential rate continuous-time}, \ic{and} the last inequality follows from the sector-bounded property of each $ \Delta^{(i)}_{\rho}$ (Lemma~\ref{lem:co-coercivity sector IQC}) and $\Gamma D = 0$.
\hiro{Then, $\| \tilde{z}_{\rho} (t)\|^2 \leq \alpha_{1}^{-1} V_{\rho} (t) \leq \alpha^{-1} V_{\rho} (0) \leq \alpha_{1}^{-1} \alpha_{2} \|\tilde{z}_{\rho}(0)\|^2$.}
Therefore, the interconnection of $G_{\rho}$ and $\Delta_{\rho}$ is linearly stable.
\lmm{As a result, we have $\| z - z^{\textup{opt}}\| \leq e^{-\rho t} c \|z_0 - z^{\textup{opt}}\|$, for some $c > 0$.} Since $z = \nabla \phi (x)$ and $\phi \in S(\mu_{\phi}, L_{\phi}) $, we can obtain $\| x - x^{\textup{opt}}\| \leq c \frac{L_{\phi}}{\mu_{\phi}} e^{-\rho t} \| x_0 - x^{\textup{opt}}\|$,  indicating that the variable $x$ also converges exponentially with rate $\rho$.
\revise{Next, to obtain the largest lower bound on the convergence rate, let us select \lr{$P = p
I_d > 0$, $p\in\mathbb{R}$}  and expand \eqref{eq:LMI exponential rate continuous-time} to obtain}
\begin{align}\label{eq: LMI with solution}
\hspace{-2mm}
\begin{bmatrix}
		2(\rho - \eta \mu_{f} \mu_{\bar{\phi}})p & * & \eta \mu_f p - q_2 - 2 \rho \gamma +\eta \gamma \mu_{f} \mu_{\bar{\phi}} \\
 \eta p - q_1 \mu_{\bar{\phi}}  & - \frac{2 q_1}{k_1} & * \\
 * & q_1 - \eta \gamma & - 2\left( \frac{q_2}{k_2} + \eta \gamma \mu_f \right)
	\end{bmatrix} \hspace{-1mm} \leq 0.
\end{align}
\hiro{
It can be verified that
\ir{ 
the above LMI always holds when}
$\rho = \eta \mu_f \mu_{\bar{\phi}}$, $p = \gamma \mu_{\bar{\phi}}$, $q_1 = \eta \gamma$, $q_2 = 0$ and $\gamma > 0$.}
Note that such $\rho$ 
\iclr{is also a tight lower bound on the convergence rate}, since it is the largest lower bound on the convergence rate \iclr{in continuous} time for \mmli{any quadratic functions} $f \in S(\mu_f, L_f)$ and $\bar{\phi} \in S(\mu_{\bar{\phi}}, L_{\bar{\phi}})$, \lir{i.e.\hu{,} one can construct quadratic functions $f, \bar{\phi}$ in these classes that achieve this rate bound,
\hiro{e.g., quadratic functions $f$, $\phi$} as in Proposition \ref{prop: MD quadratic convergence rate} with $F=\mu_f I_d$, $\Phi=L_\phi I_d$.}
}

\hiro{
\begin{remark}
\hirot{\lir{We would like to note that there is a direct}
link between the Popov criterion \lir{(Lemma~\ref{lem:popov})}  and the LMI \eqref{eq:LMI exponential rate continuous-time} in Theorem~\ref{thm:convergence rate continuous-time}.
Firstly, it} follows from \cite{hu2016exponential} that $\Delta_{\rho}$ \lir{used in the proof Theorem~\ref{thm:convergence rate continuous-time}} satisfies the Sector IQC in \eqref{eq:IQC sector bounded}, and the following modified Popov IQC
\begin{align*}
\Pi_{P, \rho} (j \omega) =
\begin{bmatrix}
			0_d & (\rho - j \omega) I _d \\ ( \rho + j\omega) I_d & 0_d
		\end{bmatrix}.
\end{align*}
\hiro{It can be established from \lir{the proof of} Theorem~\ref{thm:IQC theorem}} that the closed-loop system of $G_{\rho}(s)$ and $\Delta_{\rho}$ is stable if
\begin{align}\label{eq: modified Popov criterion}
\left\{ Q K^{-1}  - ( Q +  \rho \Gamma + j \omega \Gamma) G_{\rho}(j \omega) \right\}_{\textup{H}} \geq \delta I,  ~\forall \omega \in \mathbb{R}
\end{align}
for some $\delta >0$.
Using the KYP lemma (\cite{rantzer1996kalman}), we have $\left\{ Q K^{-1} - \hirot{\frac{\delta}{2} I } - (\mm{ Q + \rho \Gamma } + j \omega \Gamma) G_{\rho} (j\omega ) \right\}_{\hiro{\textup{H}}} \geq 0$ if and only if
\begin{align}\label{eq: strict LMI with delta}
		\begin{bmatrix}
		P\tilde{A}_{\hirot{\rho}} + \tilde{A}_{\hirot{\rho}}^TP  & P \tilde{B} - \tilde{C}_{\hirot{\rho}}^T \\
		* & - \left( \tilde{D} + \tilde{D}^T \right) \hirot{ - \delta I}
	\end{bmatrix}
	\leq 0.
\end{align}
where $\tilde{A}_{\rho} = A + \rho I$, $\tilde{B} = - B$, $\tilde{C}_{\hirot{\rho}} = \mm{(Q + 2 \rho \Gamma)} C + \Gamma C A$, and $\tilde{D} = - Q D + Q K^{-1} - \Gamma C B$.
Note that \eqref{eq: modified Popov criterion} reduces to the condition in Lemma~\ref{lem:popov} when $\rho = 0$.
\lir{It should be noted that in \eqref{eq:LMI exponential rate continuous-time} we have $\delta=0$ since linear stability of the interconnection of $G_{\rho}(s)$ and $\Delta_{\rho}$ is sufficient for the proof of Theorem \ref{thm:convergence rate continuous-time}.}
\end{remark}
}

\revise{
\section{Proof of Proposition~\ref{prop: MD quadratic convergence rate}}\label{proof of prop: MD quadratic convergence rate}
The conjugate function \hiro{of $\phi$}
can be readily derived: $\bar{\phi} (z) = \frac{1}{2} z^T \bar{\Phi} z$, where $\bar{\Phi} = \Phi^{-1}$.  Thus, $\mu_{\bar{\phi}} I_d \leq \bar{\Phi} \leq L_{\bar{\phi}} I_d$.
Following \eqref{eq:MD discrete-time composition}, the discrete-time MD algorithm is written as
\begin{align*}
z_{k+1} = z_{k} - \eta F  \bar{\Phi}  z_{k} - \eta p
\end{align*}
The error dynamics with respect to the optimal solution $z^{\text{opt}}$ is
\begin{align}\label{eq:MD algorithm quadratic functions}
\tilde{z}_{k+1} = \tilde{z}_{k} - \eta F  \bar{\Phi}  \tilde{z}_{k} = (I_d -  \eta F  \bar{\Phi} ) \tilde{z}_{k}.
\end{align}
Note that $F \bar{\Phi} = F^{1/2} \left( F^{1/2}  \bar{\Phi} F^{1/2} \right) F^{-1/2}$, where $F^{1/2}$ is the square root of $F$, meaning that $F  \bar{\Phi}$ is similar to the positive definite matrix $F^{1/2}  \bar{\Phi} F^{1/2} $ and thus all its eigenvalues are real and positive.
Let $\lambda$ be an eigenvalue of $F^{1/2}  \bar{\Phi} F^{1/2}$ associated with eigenvector $x$. Then, we have
$
x^T F^{1/2}  \bar{\Phi} F^{1/2} x = \lambda x^T x.
$
Let $q  = F^{1/2} x$,  it follows that
\begin{align*}
q^T \bar{\Phi} q  = \lambda q^T F^{-1} q.
\end{align*}
We know that $\mu_{\bar{\phi}} \|q \|^2 \leq q^T  \bar{\Phi} q \leq L_{\bar{\phi}} \| q \|^2$, and $1/L_{f} \|q\|^2 \leq q^T F^{-1} q \leq 1/m_{f} \|q\|^2$.  Substituting them back to the previous equation, we can obtain the bound for the spectrum of $F \bar{\Phi}$:
\begin{align}\label{eq: spectrum bound}
m_{f} m_{\bar{\phi}} \leq \lambda \leq L_{f} L_{\bar{\phi}}.
\end{align}
As a result, there exists a linear transformation $\zeta = Q\tilde{z}$, such that \eqref{eq:MD algorithm quadratic functions} is transformed to
\begin{align}\label{eq: MD algorithm quadratic transformed}
{\zeta}_{k+1} = (I_d - \eta \Lambda) \zeta_{k}
\end{align}
where $\Lambda$ is a diagonal positive matrix whose elements are eigenvalues of $F \bar{\Phi} $, bounded by \eqref{eq: spectrum bound}.
Since system \eqref{eq: MD algorithm quadratic transformed} shares the same form as the GD method applied to the quadratic convex function $\frac{1}{2} \zeta^T \Lambda \zeta$, the optimal stepsize and a tight convergence rate \iclr{bound} can be directly inferred from the analysis of GD method \cite[Theorem 2.1.14]{nesterov2003introductory}.
Therefore, the optimal stepsize is given by $\eta =  \frac{2}{L_f L_{\bar{\phi}} + \mu_f \mu_{\bar{\phi}}}$ and the tight convergence rate \iclr{bound} is $\rho = \frac{ \kappa_f \kappa_{\bar{\phi}} -1 }{ \kappa_f \kappa_{\bar{\phi}} + 1}$ for both \eqref{eq: MD algorithm quadratic transformed} and its linear transformation \eqref{eq:MD algorithm quadratic functions}.

\section{Matrices $\hat{A}$, $\hat{B}$, $\hat{C}$ and $\hat{D}$ in \eqref{eq:LMI IQC discrete-time constrained MD}}\label{Matrices for LMI}
\begin{figure*}[htbp]
\hiro{\begin{align}\label{eq:LMI matrices constrained MD}
	\hat{A} = &
	\begin{bmatrix}
		A \\
		B_{{\Psi_{w,f}}}^{y} C_1 & A_{\Psi_{w,f}}\\
	 	B_{{\Psi_{w, \bar{\phi}}}}^{y} C_2 &  & A_{\Psi_{w,\bar{\phi}}}\\
	 	B_{{\Psi_{w, T}}}^{y} C_3 & & & A_{\Psi_{w,T}} \\
	 	B_{{\Psi_{w, \bar{\phi}}}}^{y} C_4 & & & & A_{\Psi_{w, \bar{\phi} }}\\
	 	B_{{\Psi_{w, \bar{\phi}}}}^{y} (C_2 - C_4) & & & & & A_{\Psi_{w, \bar{\phi}}}
	\end{bmatrix},
	\\
	\hat{B} = &
	\begin{bmatrix}
    B_1 & B_2 & B_3 & B_4\\
    B_{\Psi_{w, f}}^{y} D_{11} + B_{\Psi_{w, f}}^{u} & B_{\Psi_{w, f}}^{y} D_{12} & B_{\Psi_{w, f}}^{y} D_{13} & B_{\Psi_{w, f}}^{y} D_{14}\\
    B_{\Psi_{w, \bar{\phi}}}^{y} D_{21} & B_{\Psi_{w, \bar{\phi}}}^{y} D_{22} + B_{\Psi_{w, \bar{\phi}}}^{u} & B_{\Psi_{w, \bar{\phi}}}^{y} D_{23} & B_{\Psi_{w, \bar{\phi}}}^{y} D_{24}\\
    B_{\Psi_{w, T}}^{y} D_{31} & B_{\Psi_{w, T}}^{y} D_{32} & B_{\Psi_{w, T}}^{y} D_{33} + B_{\Psi_{w, T}}^{u} & B_{\Psi_{w, T}}^{y} D_{34}\\
    B_{\Psi_{w, \bar{\phi}}}^{y} D_{41} & B_{\Psi_{w, \bar{\phi}}}^{y} D_{42}  & B_{\Psi_{w, \bar{\phi}}}^{y} D_{43} & B_{\Psi_{w, \bar{\phi}}}^{y} D_{44} + B_{\Psi_{w, \bar{\phi}}}^{u}\\
    B_{\Psi_{w, \bar{\phi}}}^{y} (D_{21} - D_{41}) & B_{\Psi_{w, \bar{\phi}}}^{y} (D_{22} - D_{42}) + B_{\Psi_{w, \bar{\phi}}}^{u} & B_{\Psi_{w, \bar{\phi}}}^{y} (D_{23} - D_{43}) & B_{\Psi_{w, \bar{\phi}}}^{y} (D_{24} - D_{44}) - B_{\Psi_{w, \bar{\phi}}}^{u}
	\end{bmatrix},
	 \\
\hat{C} = &
	\begin{bmatrix}
		D_{\Psi_{s,f}}^{y} C_1 \\
		D_{\Psi_{w,f}}^{y} C_1 & C_{\Psi_{w,f}} \\
		D_{\Psi_{s, \bar{\phi}}}^{y} C_2 \\
		D_{\Psi_{w, \bar{\phi}}}^{y} C_2 & & C_{\Psi_{w, \bar{\phi}}}\\
		D_{\Psi_{s, T}}^{y} C_3 \\
		D_{\Psi_{w, T}}^{y} C_3 & & & C_{\Psi_{w, T}}\\
		D_{\Psi_{s, \bar{\phi}}}^{y} C_4 \\
		D_{\Psi_{w, \bar{\phi}}}^{y} C_4 & & & & C_{\Psi_{w, \bar{\phi}}}\\
		D_{\Psi_{s, \bar{\phi}}}^{y} (C_2 - C_4) \\
		D_{\Psi_{w, \bar{\phi}}}^{y} (C_2 - C_4) & & & & & C_{\Psi_{w, \bar{\phi}}}
	\end{bmatrix},
	\\
	\hat{D} = &
	\begin{bmatrix}
		D_{\Psi_{s, f}}^{y} D_{11} + D_{\Psi_{s, f}}^{u} & D_{\Psi_{s, f}}^{y} D_{12} & D_{\Psi_{s, f}}^{y} D_{13} & D_{\Psi_{s, f}}^{y} D_{14}\\
		D_{\Psi_{w, f}}^{y} D_{11} + D_{\Psi_{w, f}}^{u} & D_{\Psi_{w, f}}^{y} D_{12} & D_{\Psi_{w, f}}^{y} D_{13} & D_{\Psi_{w, f}}^{y} D_{14}\\
		D_{\Psi_{s, \bar{\phi}}}^{y} D_{21} & D_{\Psi_{s, \bar{\phi}}}^{y} D_{22} + D_{\Psi_{s, \bar{\phi}}}^{u} & D_{\Psi_{s, \bar{\phi}}}^{y} D_{23} & D_{\Psi_{s, \bar{\phi}}}^{y} D_{24}\\
		D_{\Psi_{w, \bar{\phi}}}^{y} D_{21} & D_{\Psi_{w, \bar{\phi}}}^{y} D_{22} + D_{\Psi_{w, \bar{\phi}}}^{u} & D_{\Psi_{w, \bar{\phi}}}^{y} D_{23} & D_{\Psi_{w, \bar{\phi}}}^{y} D_{24} \\
		D_{\Psi_{s, T}}^{y} D_{31} & D_{\Psi_{s, T }}^{y} D_{32} + D_{\Psi_{s, T}}^{u} & D_{\Psi_{s, T}}^{y} D_{33} + D_{\Psi_{s, T}}^{u}  & D_{\Psi_{s, T}}^{y} D_{34}\\
		D_{\Psi_{w, T}}^{y} D_{31} & D_{\Psi_{w, T }}^{y} D_{32} + D_{\Psi_{w, T}}^{u} & D_{\Psi_{w, T}}^{y} D_{33} + D_{\Psi_{w, T}}^{u}  & D_{\Psi_{w, T}}^{y} D_{34}\\
		D_{\Psi_{s, \bar{\phi}}}^{y} D_{41} & D_{\Psi_{s, \bar{\phi}}}^{y} D_{42} & D_{\Psi_{s, \bar{\phi}}}^{y} D_{43} & D_{\Psi_{s, \bar{\phi}}}^{y} D_{44} + D_{\Psi_{s, \bar{\phi}}}^{u}\\
		D_{\Psi_{w, \bar{\phi}}}^{y} D_{41} & D_{\Psi_{w, \bar{\phi}}}^{y} D_{42} & D_{\Psi_{w, \bar{\phi}}}^{y} D_{43} & D_{\Psi_{w, \bar{\phi}}}^{y} D_{44} + D_{\Psi_{w, \bar{\phi}}}^{u}\\
		D_{\Psi_{s, \bar{\phi}}}^{y} (D_{21} - D_{41}) & D_{\Psi_{s, \bar{\phi}}}^{y} (D_{22} - D_{42}) + D_{\Psi_{s, \bar{\phi}}}^{u} & D_{\Psi_{s, \bar{\phi}}}^{y} (D_{23} - D_{43}) & D_{\Psi_{s, \bar{\phi}}}^{y} (D_{24} - D_{44}) - D_{\Psi_{s, \bar{\phi}}}^{u}\\
		D_{\Psi_{w, \bar{\phi}}}^{y} (D_{21} - D_{41}) & D_{\Psi_{w, \bar{\phi}}}^{y} (D_{22} - D_{42}) + D_{\Psi_{w, \bar{\phi}}}^{u} & D_{\Psi_{w, \bar{\phi}}}^{y} (D_{23} - D_{43}) & D_{\Psi_{w, \bar{\phi}}}^{y} (D_{24} - D_{44}) - D_{\Psi_{w, \bar{\phi}}}^{u}
	\end{bmatrix}.
\end{align}}
\end{figure*}

\end{document}